\documentclass[11pt]{amsart}
\usepackage{euscript}
\usepackage{amssymb}
\usepackage{amsmath}
\usepackage{epic}
\usepackage[matrix,arrow,curve]{xy}

\newtheoremstyle{my}{1.5em}{0.5em}{\em}{}{\sc}{.}{0.5em}{}

\theoremstyle{my}
\newtheorem{thm}{Theorem}
\newtheorem{Theorem}[thm]{Theorem}
\newtheorem*{Theorem*}{Theorem}

\newtheorem*{corollary*}{Corollary}
\newtheorem{Lemma}[thm]{Lemma}

\newtheorem*{conjecture*}{Conjecture}

\newtheorem*{question*}{Question}

\newtheorem*{definitions*}{Definitions}

\newtheorem*{rem*}{Remark}
\newtheorem{Remark}[thm]{Remark}
\newtheorem*{remark*}{Remark}

\newtheorem*{remarks*}{Remarks}
\newtheorem*{example*}{Example}
\newtheorem{Addendum}[thm]{Addendum}

\newtheorem*{examples*}{Examples}

\newcommand{\ignore}[1]{...}

\newcommand{\R}{\mathbb{R}}
\newcommand{\Z}{\mathbb{Z}}

\newcommand{\C}{\mathbb{C}}

\newcommand{\iso}{\cong}           
\newcommand{\htp}{\simeq}          

\newcommand{\CP}[1]{\C {\mathrm P}^{#1}}
\newcommand{\RP}[1]{\R {\mathrm P}^{#1}}
\newcommand{\leftsc}{\langle}
\newcommand{\rightsc}{\rangle}

\newcommand{\gapprox}{\gtrsim}
\newcommand{\lapprox}{\lesssim}

\newcommand{\re}{\mathrm{re}}
\newcommand{\im}{\mathrm{im}}

\renewcommand{\o}{\omega}
\renewcommand{\O}{\mathcal{O}}

\newcommand{\K}{\mathbb{K}}
\renewcommand{\index}{\mathrm{index}\,}

\newcommand{\A}{\mathcal{A}}
\newcommand{\MM}{\mathcal{M}}
\newcommand{\NN}{\mathcal{N}}
\newcommand{\PP}{\mathcal{P}}
\newcommand{\QQ}{\mathcal{Q}}
\renewcommand{\SS}{\mathcal{S}}
\newcommand{\F}{\mathcal{F}}

\newcommand{\CC}{\mathcal{C}}

\title[Exact Lagrangian submanifolds]{Exact Lagrangian submanifolds in
  simply-connected cotangent bundles}
\author{Kenji Fukaya, \, Paul Seidel, \, Ivan Smith}
\date{January 2007}
\begin{document}
\maketitle
\begin{abstract}
We consider exact Lagrangian submanifolds in cotangent bundles. Under
certain additional restrictions (triviality of the fundamental group
of the cotangent bundle, and of the Maslov class and second
Stiefel-Whitney class of the Lagrangian submanifold) we prove such
submanifolds are  Floer-cohomologically indistinguishable
from the zero-section.  This implies strong restrictions on their
topology. An essentially equivalent result was recently proved
independently by Nadler \cite{nadler06}, using a different approach.
\end{abstract}

\section{Introduction\label{sec:1}}

This paper is concerned with the topology of Lagrangian submanifolds in cotangent bundles. Take a closed manifold $N$ (throughout the entire paper, the convention is that all manifolds are assumed to be connected). Equip the cotangent bundle $T^*N$ with the standard symplectic structure. We will be interested in closed exact Lagrangian submanifolds in $T^*N$.

\begin{Theorem} \label{th:main}
Suppose that $N$ is simply-connected and spin. Fix a coefficient field $\K$ of characteristic $\neq 2$. Let $L \subset T^*N$ be a closed exact Lagrangian submanifold, which is spin and whose Maslov class $m_L \in H^1(L)$ vanishes. Then the projection $L \rightarrow N$ has degree $\pm 1$, and induces an isomorphism $H^*(N;\K) \rightarrow H^*(L;\K)$. Moreover, if $L_0,L_1$ are two submanifolds satisfying the same conditions, and intersecting transversally, then $|L_0 \cap L_1| \geq \dim\, H^*(N;\K)$, where the right hand side is the sum of the Betti numbers with $\K$-coefficients.
\end{Theorem}

A well-known conjecture says that all closed exact $L \subset T^*N$ should be isotopic to the zero-section (where the isotopy goes through exact Lagrangian submanifolds). In this form, the conjecture (sometimes called the ``nearby Lagrangian problem'') seems to be beyond the reach of present technology, but there is a long history of partial results. Surjectivity of the projection $N \rightarrow L$ was proved in one of the first papers on the subject
\cite{lalonde-sikorav91}. Furthermore, many non-embedding results are
known for special classes of manifolds $N$ or $L$; besides the
reference already quoted, see \cite{viterbo94, viterbo97a, viterbo97b,
  buhovsky03, hind03, seidel04c} (even this is a non-exhaustive
list). Together, these papers use a wide variety of approaches, and
the consequent statements vary considerably in strength (sometimes far
outstripping what one can get by applying Theorem \ref{th:main}; for
instance, results in \cite{viterbo94} and \cite{hind03} prove that
every oriented exact $L \subset T^*S^2$ is indeed Lagrangian isotopic
to the zero-section). This diversity is one of the aspects making this
an interesting problem to study.

To bring the story to a close, there is an important very recent paper
of Nadler \cite{nadler06} in which, building on work of Nadler-Zaslow
\cite{nadler-zaslow06} and  Fukaya-Oh
\cite{fukaya-oh98}, he obtains a result essentially equivalent to our
Theorem \ref{th:main}. Nadler's argument is somewhat different from
the one used here (there is also yet another approach, due to the
authors of the present paper, which remains so far unpublished). We
emphasize that Nadler's work and ours were carried out entirely
independently of each other. Nevertheless, there are many similarities
on a philosophical level; notably, the use of the Fukaya category of
$T^*N$, enlarged by admitting certain non-compact Lagrangian
submanifolds, and of decompositions of the diagonal.  We postpone a
more detailed comparison (and a discussion of the situation when $\pi_1(N)\neq 0$) to \cite{fukaya-seidel-smith07b}.

This paper is organized as follows. The rest of Section \ref{sec:1}
gives a complete account of the proof of Theorem \ref{th:main},
assuming certain auxiliary theorems which are then addressed in the
subsequent sections. Mostly, the proofs of those auxiliary results are
quite self-contained, relying only on classical Floer homology theory
and algebraic geometry. However, there is one notable exception,
Theorem \ref{th:ss}, which belongs to the general theory of Lefschetz
fibrations. The use of Lefschetz fibrations is not, of course,
intrinsic to the nearby Lagrangian problem. Indeed, a large part of
the paper is devoted to going from one framework (cotangent bundles)
to the other (affine algebraic varieties) and back. Nevertheless, this
turns out to be a price worth paying, because the Fukaya
categories associated to Lefschetz fibrations belong to a particularly
benign class (directed $A_\infty$-categories, up to derived
equivalence).

{\em Acknowledgements}. The first author was partially supported by JSPS Grant-in-Aid for Scientific Research 18104001. The second author was partially supported by
NSF grant DMS-0405516. The third author was partially supported by
EPSRC grant EP/C535995/1.

\subsection{Lefschetz fibrations on affine varieties\label{subsec:lefschetz}}

Let $\bar{X}$ be a smooth $n$-dimensional complex projective variety, equipped with an ample line bundle $\O_{\bar{X}}(1)$. Let $t \in H^0(\O_{\bar{X}}(1))$ be a holomorphic section, which defines a normal crossing divisor (possibly with multiplicities) $Y = t^{-1}(0)$. The complement $X = \bar{X} \setminus Y$ is an affine variety. We fix a hermitian metric on $\O_{\bar{X}}(1)$ whose curvature is positive, hence defines a K{\"a}hler form. The restriction of this form to $X$ can be written as $\omega = -dd^c h$, where $h = -\log ||t||^2$ is the K{\"a}hler potential. In particular, the symplectic form has a canonical primitive, namely $\theta = -d^ch$.

Take another section $s \in H^0(\O_{\bar{X}}(1))$, which is not a multiple of $t$. We say that the function $p = s/t: X \rightarrow \C$ is a {\em Lefschetz fibration} if the following conditions are satisfied:
\begin{itemize} \itemsep0.5em
\item
$p$ has only nondegenerate critical points, and there is at most one such point in every fibre.
\item
$s^{-1}(0)$ is (reduced and) smooth near $Y$, and intersects each stratum of $Y$ transversally.
\end{itemize}
The first condition is the standard Lefschetz property, and the second one ensures that the fibres are well-behaved at infinity. More precisely, if we consider the symplectic connection defined by $\o$, then:

\begin{Lemma} \label{th:parallel}
$p: X \rightarrow \C$ has well-defined symplectic parallel transport maps (away from the singular fibres).
\end{Lemma}

We need to review some more terminology from Picard-Lefschetz theory. A {\em vanishing path} is an embedded path $\gamma: [0;\infty) \rightarrow \C$ such that $\gamma(0)$ is a critical value of $p$, all the other $\gamma(r)$ are regular values, and $\gamma'(r)$ is constant for $r \gg 0$ (which means that the path eventually becomes a straight half-line). For any such path there is an associated {\em Lefschetz thimble} $\Delta_{\gamma} \subset X$, which is an open Lagrangian disc projecting properly to $\gamma$. In fact, $\gamma^{-1} \circ p|\Delta_{\gamma}$ is the standard exhausting Morse function on the open disc, with a single nondegenerate minimum that lies precisely at the unique singular point of $p^{-1}(\gamma(0))$. We refer to \cite{seidel01} for details.

Let's label the critical values by $\{z_1,\dots,z_m\}$. A pair of {\em dual bases of vanishing paths} consists of collections $\{\gamma_j\}$, $\{\gamma^!_j\}$ indexed by $1 \leq j \leq m$, with the following properties:
\begin{itemize} \itemsep0.5em
\item $\gamma_j(0) = \gamma_j^!(0) = z_j$ for all $j$, with transverse intersection at that point, and no other intersection points; except for this, any two of the paths $\{\gamma_1,\dots,\gamma_m, \gamma_1^!,\dots,\gamma_m^!\}$ are disjoint.
\item For $r \gg 0$, $\gamma_j(r) = c_j - ir$ for some constants $c_j \in \C$, which are such that their real parts $\re\,c_j$ are increasing with $j$. Similarly, $\gamma_j^!(r) = c_j^! + ir$, with $\re\, c_j^!$ increasing with $j$.
\end{itemize}
These bases give rise to two collections of Lefschetz thimbles, which we usually denote by $\{\Delta_j\}$ and $\{\Delta_j^!\}$. By construction, $\Delta_j$ intersects $\Delta_j^!$ precisely at the unique critical point $x_j \in p^{-1}(z_j)$ (and the intersection is transverse there); otherwise, the Lefschetz thimbles are all mutually disjoint. Note that $\Delta_1$ is always isotopic to $\Delta_1^!$, and similarly $\Delta_m$ to $\Delta_m^!$, just because the relevant paths can be moved into each other.

\subsection{Floer cohomology\label{subsec:maslov}}

Take an affine variety $X = \bar{X} \setminus Y$ of the kind considered above, with its symplectic form $\o = d\theta$. In addition, assume that $\bar{X}$ carries a meromorphic complex $n$-form whose zeros and poles lie entirely inside $Y$. Restriction of that form to $X$ then yields a holomorphic volume form, which we denote by $\eta$. For any Lagrangian submanifold $L \subset X$, we then have a class $[\theta|L] \in H^1(L;\R)$, as well as the Maslov class $m_L \in H^1(L)$. We define {\em admissible Lagrangian submanifolds} to be those $L$ which are:
\begin{itemize}
\itemsep0.5em
\item exact, meaning that $[\theta|L] = 0$;
\item spin, and in fact come with a choice of spin structure; and
\item have zero Maslov class.
\end{itemize}
The last property allows us to choose a grading of $L$, in the terminology of \cite{seidel99}. We will denote the resulting graded Lagrangian submanifold by $\tilde{L}$. For any pair of closed submanifolds of this kind, one has a well-defined Floer cohomology group $HF^*(\tilde{L}_0,\tilde{L}_1)$, which is a $\Z$-graded vector space (over a coefficient field $\K$, which can be chosen arbitrarily). One can also allow certain non-compact Lagrangian submanifolds. Namely, if $X$ carries a Lefschetz fibration in the sense defined above, then one or both $\tilde{L}_k$ may be Lefschetz thimbles (these are contractible, hence automatically admissible), except for one condition: in the case of two Lefschetz thimbles, we require that the associated paths should only intersect in a compact subset of $\C$. The resulting Floer cohomology groups have the usual properties: for $\tilde{L}_0 = \tilde{L}_1 = \tilde{L}$, one has $HF^*(\tilde{L},\tilde{L}) \iso H^*(L;\K)$; isotopy invariance holds within the class of Lagrangian submanifolds which are allowed; there is a Poincar{\'e} duality isomorphism $HF^*(\tilde{L}_1,\tilde{L}_0) \iso HF^{n-*}(\tilde{L}_0,\tilde{L}_1)^\vee$; and also an associative product $HF^*(\tilde{L}_1,\tilde{L}_2) \otimes HF^*(\tilde{L}_0,\tilde{L}_1) \rightarrow HF^*(\tilde{L}_0,\tilde{L}_2)$.

\subsection{A spectral sequence}
Fix dual bases of vanishing paths, and consider the associated Lefschetz thimbles. Recall that whenever two graded Lagrangian submanifolds $\tilde{L}_0,\tilde{L}_1$ intersect transversally at a point $x$, one can define the {\em absolute Maslov index}
\begin{equation} \label{eq:absolute}
i(\tilde{L}_0,\tilde{L}_1;x) \in \Z,
\end{equation}
which determines the degree in which this point will contribute to the Floer cochain group $CF^*(\tilde{L}_0,\tilde{L}_1)$. In our case, our convention is to choose gradings so that the unique intersection point $x_j \in p^{-1}(z_j)$ satisfies
\begin{equation} \label{eq:dual-convention}
i(\tilde\Delta_j^!,\tilde\Delta_j;x_j) = 0.
\end{equation}
This implies that $HF^*(\tilde\Delta_j^!,\tilde\Delta_j)$ is $\K$ in degree $0$, and trivial in other degrees. There is some residual ambiguity, since one can change the grading of both vanishing cycles by the same amount, but that will be irrelevant for our purpose. 

\begin{Lemma}\label{lem:gradingconvention}
This grading convention is compatible with the isotopy $\Delta_m \htp \Delta_m^!$ (but not with $\Delta_1 \htp \Delta_1^!$).
\end{Lemma}

\begin{Theorem} \label{th:ss}
Suppose that $char(\K) \neq 2$. Let $(\tilde{L}_0,\tilde{L}_1)$ be two closed admissible Lagrangian submanifolds in $X$, equipped with gradings. Then there is a spectral sequence converging to $HF^*(\tilde{L}_0,\tilde{L}_1)$, whose starting page is
\begin{equation} \label{eq:e1}
E_1^{jk} =
\Big(HF^*(\tilde\Delta_j^!,\tilde{L}_1) \otimes HF^*(\tilde{L}_0,\tilde\Delta_j) \Big)^{j+k}.
\end{equation}
\end{Theorem}

It is implicit in the statement that the columns $E_1^{j\ast}$ with $j \leq 0$ or $j > m$ are zero. The result is quoted from \cite[Corollary 18.27]{seidel04} with some notational changes, which we will comment on later. More importantly, we will need two additional properties of this spectral sequence. First, like any spectral sequence with bounded starting term, ours comes with a right-sided {\em edge homomorphism}, which is a map $E_1^{m,*-m} \rightarrow HF^*(\tilde{L}_0,\tilde{L}_1)$.

\begin{Addendum} \label{th:edge}
Up to a nonzero multiplicative constant, the edge homomorphism is the composition of the isomorphism $HF^*(\tilde\Delta_m,\tilde{L}_1) \iso HF^*(\tilde\Delta_m^!,\tilde{L}_1)$ obtained from the isotopy $\Delta_m \htp \Delta_m^!$, and the multiplication $HF^*(\tilde\Delta_{m},\tilde{L}_1) \otimes HF^*(\tilde{L}_0,\tilde\Delta_m) \longrightarrow HF^*(\tilde{L}_0,\tilde{L}_1)$.
\end{Addendum}

The second observation can be motivated as follows. Suppose for a second that $HF^*(\tilde{L}_0,\tilde\Delta_j) = 0$ for all $j$. Then, the spectral sequence would imply that $HF^*(\tilde{L}_0,\tilde{L}_1) = 0$ for all $\tilde{L}_1$. This can never actually happen, since $HF^*(\tilde{L}_0,\tilde{L}_0) \neq 0$, but there is a useful relative version (a kind of Whitehead theorem):

\begin{Addendum} \label{th:cone}
Suppose that there is a $c \in HF^0(\tilde{L}_1,\tilde{L}_0)$ such that the product with $c$ is an isomorphism $HF^*(\tilde{L}_0,\tilde\Delta_j) \rightarrow HF^*(\tilde{L}_1,\tilde\Delta_j)$ for all $j$. Then, the same product is an isomorphism $HF^*(\tilde{L}_0,\tilde{L}) \rightarrow HF^*(\tilde{L}_1,\tilde{L})$ for all closed admissible $L$.
\end{Addendum}

\subsection{Real Lefschetz fibrations\label{subsec:squish}}
We will be looking at {\em Lefschetz fibrations with real structures}. This means that
\begin{itemize} \itemsep0.5em
\item
$\bar{X}$ carries an anti-holomorphic involution, which comes with a lift to the line bundle $\O_{\bar{X}}(1)$, preserving the hermitian metric. The section $t$ should be defined over $\R$, which means that it is invariant under the involution;
\item
if we restrict the involution to $X$, then its fixed part $X_\R$ should be compact (and connected);
\item
the section $s$ defining the Lefschetz fibration should also be defined over $\R$.
\end{itemize}
Note that $X_\R$ is automatically an exact Lagrangian submanifold, since both $\o$ and $\theta$ vanish on it. Restriction of $p$ to the real part yields a function $p_\R: X_\R \rightarrow \R$, which is automatically Morse. Moreover, because we do not allow more than one critical point in a given fibre, the critical points of $p_\R$ are precisely those critical points of $p$ which have real values. In addition to the properties above, we will often require the existence of a holomorphic volume form $\eta$ as before, which should again be defined over $\R$ (this means that the pullback of $\eta$ by the anti-holomorphic involution is $\bar\eta$). Then, the restriction $\eta_\R = \eta|X_\R$ is a volume form on $X_\R$.

In the presence of a real structure, it is convenient to use dual bases of vanishing paths of a particular kind. Suppose that there are $s$ real and $r$ non-real critical values, so $m = r+s$. We choose an arbitrary ordering $\{z_1,\dots,z_r\}$ of the non-real ones, and then add the real ones in their natural order as $\{z_{r+1} < \cdots < z_m\}$. Let $\{\gamma_j\}$, $\{\gamma_j^!\}$ be dual bases of vanishing paths, with $\gamma_j(0) = \gamma_j^!(0) = z_j$. We say that they are {\em compatible with the real structure} if:
\begin{itemize} \itemsep0.5em
\item for each $j \leq r$, $\gamma_j$ is disjoint from $p(X_\R)$;
\item for each $j > r$, both $\gamma_j$ and $\gamma_j^!$ intersect $p(X_\R)$ only at their starting point. Moreover, $\im\, \gamma_j'(0) < 0$ and $\im\,(\gamma_j^!)'(0) > 0$.
\end{itemize}
Figure \ref{fig:real-basis} shows how to find such bases, first in a particularly simple situation, and then in a more realistic one (the two cases are related by a diffeomorphism of $\C$ preserving $p(X_\R)$, hence are not really substantially different).
\begin{figure}[ht]
\begin{center}
\setlength{\unitlength}{0.00062500in}
\begingroup\makeatletter\ifx\SetFigFont\undefined%
\gdef\SetFigFont#1#2#3#4#5{%
  \reset@font\fontsize{#1}{#2pt}%
  \fontfamily{#3}\fontseries{#4}\fontshape{#5}%
  \selectfont}%
\fi\endgroup%
{\renewcommand{\dashlinestretch}{30}
\begin{picture}(3387,5889)(0,-10)
\put(900,762){\makebox(0,0)[lb]{{\SetFigFont{10}{12}{\rmdefault}{\mddefault}{\updefault}$\gamma_2^!$}}}
\put(3150,4512){\circle*{100}}
\put(750,5187){\circle*{100}}
\put(1275,1212){\circle*{100}}
\put(2475,1212){\circle*{100}}
\put(750,4137){\circle*{100}}
\put(1500,612){\circle*{100}}
\put(1500,1812){\circle*{100}}
\drawline(3150,5637)(3150,3687)
\dashline{60.000}(1950,5862)(1950,5562)
\dashline{60.000}(1950,3687)(1950,3387)
\dashline{60.000}(3150,3687)(3150,3387)
\dashline{60.000}(3150,5862)(3150,5562)
\dashline{60.000}(1350,5862)(1350,5562)
\dashline{60.000}(750,5862)(750,5562)
\dashline{60.000}(750,3687)(750,3387)
\dashline{60.000}(2550,312)(2550,12)
\dashline{60.000}(3075,312)(3075,12)
\dashline{60.000}(2775,2562)(2775,2262)
\dashline{60.000}(3375,2562)(3375,2262)
\drawline(1950,5637)(1950,3687)
\dottedline{45}(1950,4512)(3150,4512)
\dottedline{45}(1275,1212)(2475,1212)
\drawline(750,4137)(750,3687)
\drawline(750,5637)(750,5187)
\drawline(150,3987)(150,3687)
\dashline{60.000}(150,3687)(150,3387)
\dashline{60.000}(2175,2562)(2175,2262)
\dashline{60.000}(1500,312)(1500,12)
\drawline(1500,612)(1500,312)
\dashline{60.000}(450,312)(450,12)
\dashline{60.000}(1500,2562)(1500,2262)
\drawline(1500,1812)(1500,2337)
\drawline(1275,1212)(1278,1211)(1283,1209)
	(1293,1206)(1309,1200)(1331,1193)
	(1358,1184)(1391,1173)(1429,1160)
	(1470,1145)(1514,1130)(1559,1115)
	(1605,1099)(1650,1083)(1694,1068)
	(1736,1053)(1775,1039)(1812,1026)
	(1847,1014)(1880,1002)(1910,990)
	(1938,980)(1965,970)(1989,960)
	(2012,951)(2034,942)(2055,933)
	(2075,924)(2102,913)(2128,901)
	(2153,889)(2177,876)(2200,864)
	(2223,851)(2245,839)(2266,826)
	(2287,813)(2306,799)(2324,786)
	(2341,773)(2357,760)(2372,748)
	(2385,735)(2397,723)(2409,710)
	(2419,698)(2429,686)(2438,674)
	(2448,659)(2458,644)(2467,627)
	(2476,611)(2484,593)(2492,576)
	(2499,558)(2506,540)(2512,523)
	(2517,506)(2522,490)(2526,475)
	(2529,461)(2533,448)(2535,436)
	(2538,424)(2541,408)(2543,394)
	(2545,380)(2547,366)(2548,351)
	(2549,337)(2549,325)(2550,316)
	(2550,313)(2550,312)
\drawline(1275,1212)(1278,1213)(1284,1215)
	(1296,1219)(1314,1225)(1338,1233)
	(1368,1243)(1402,1255)(1440,1268)
	(1480,1282)(1520,1296)(1561,1310)
	(1600,1324)(1637,1337)(1672,1349)
	(1705,1361)(1735,1372)(1763,1383)
	(1789,1392)(1814,1402)(1837,1411)
	(1859,1420)(1880,1428)(1900,1437)
	(1922,1446)(1943,1456)(1964,1466)
	(1985,1475)(2006,1486)(2027,1496)
	(2048,1507)(2069,1518)(2090,1529)
	(2110,1541)(2131,1553)(2151,1565)
	(2170,1577)(2189,1589)(2208,1602)
	(2226,1614)(2244,1626)(2261,1638)
	(2277,1650)(2294,1662)(2309,1675)
	(2325,1687)(2342,1701)(2360,1715)
	(2377,1730)(2395,1746)(2413,1762)
	(2432,1778)(2450,1795)(2469,1813)
	(2487,1830)(2506,1848)(2523,1865)
	(2541,1883)(2557,1899)(2573,1916)
	(2588,1931)(2602,1946)(2616,1961)
	(2628,1974)(2639,1987)(2650,1999)
	(2664,2016)(2677,2032)(2689,2048)
	(2701,2063)(2711,2078)(2720,2093)
	(2729,2107)(2736,2120)(2743,2133)
	(2748,2145)(2753,2156)(2757,2167)
	(2760,2177)(2763,2187)(2765,2199)
	(2767,2210)(2769,2223)(2771,2238)
	(2772,2254)(2773,2272)(2774,2290)
	(2774,2308)(2775,2323)(2775,2332)
	(2775,2336)(2775,2337)
\drawline(2475,1212)(2477,1210)(2482,1206)
	(2491,1198)(2505,1186)(2523,1170)
	(2546,1150)(2572,1128)(2600,1103)
	(2630,1076)(2661,1049)(2691,1022)
	(2720,996)(2747,972)(2773,949)
	(2796,927)(2818,907)(2838,889)
	(2856,871)(2872,855)(2887,840)
	(2900,826)(2913,813)(2925,799)
	(2941,780)(2957,762)(2971,744)
	(2984,726)(2996,707)(3007,690)
	(3017,672)(3025,655)(3033,638)
	(3040,622)(3046,606)(3050,591)
	(3054,577)(3057,563)(3060,550)
	(3063,537)(3065,522)(3067,507)
	(3068,491)(3070,474)(3071,454)
	(3072,433)(3073,410)(3074,386)
	(3074,363)(3075,342)(3075,326)
	(3075,317)(3075,313)(3075,312)
\drawline(2475,1212)(2477,1214)(2482,1218)
	(2492,1226)(2505,1238)(2524,1254)
	(2547,1274)(2574,1298)(2603,1323)
	(2634,1351)(2666,1378)(2697,1406)
	(2728,1433)(2757,1459)(2784,1483)
	(2810,1506)(2833,1527)(2855,1547)
	(2876,1566)(2895,1584)(2913,1601)
	(2930,1618)(2947,1634)(2963,1649)
	(2981,1668)(3000,1687)(3018,1706)
	(3036,1725)(3054,1745)(3071,1765)
	(3089,1785)(3106,1805)(3123,1825)
	(3140,1844)(3155,1864)(3171,1883)
	(3185,1901)(3198,1919)(3211,1936)
	(3223,1953)(3234,1969)(3244,1984)
	(3254,1998)(3263,2012)(3274,2031)
	(3285,2050)(3296,2069)(3305,2088)
	(3315,2107)(3323,2125)(3331,2143)
	(3337,2159)(3343,2175)(3349,2190)
	(3353,2203)(3357,2216)(3360,2227)
	(3363,2237)(3366,2250)(3368,2262)
	(3370,2274)(3372,2287)(3373,2300)
	(3374,2313)(3374,2325)(3375,2333)
	(3375,2336)(3375,2337)
\drawline(750,4137)(750,4140)(751,4148)
	(752,4160)(754,4178)(757,4199)
	(761,4223)(765,4248)(770,4273)
	(776,4297)(783,4320)(791,4343)
	(800,4365)(811,4388)(823,4412)
	(837,4437)(848,4455)(859,4473)
	(871,4492)(884,4512)(898,4533)
	(913,4555)(928,4578)(945,4601)
	(961,4625)(979,4650)(996,4674)
	(1014,4700)(1032,4725)(1050,4749)
	(1068,4774)(1085,4798)(1102,4821)
	(1118,4844)(1133,4866)(1148,4887)
	(1162,4907)(1176,4926)(1188,4944)
	(1200,4962)(1215,4984)(1228,5006)
	(1241,5027)(1253,5047)(1264,5067)
	(1274,5086)(1284,5105)(1293,5123)
	(1300,5141)(1307,5158)(1313,5175)
	(1319,5190)(1323,5206)(1327,5220)
	(1330,5234)(1333,5248)(1335,5261)
	(1338,5274)(1340,5290)(1341,5306)
	(1343,5322)(1344,5340)(1345,5361)
	(1346,5383)(1347,5407)(1348,5433)
	(1348,5461)(1349,5488)(1349,5513)
	(1350,5534)(1350,5549)(1350,5558)
	(1350,5561)(1350,5562)
\drawline(150,3987)(150,3990)(151,3998)
	(152,4010)(154,4028)(157,4049)
	(161,4073)(165,4098)(170,4123)
	(176,4147)(183,4170)(191,4193)
	(200,4215)(211,4238)(223,4262)
	(237,4287)(248,4305)(259,4323)
	(271,4342)(284,4362)(298,4383)
	(313,4405)(328,4428)(345,4451)
	(361,4475)(379,4500)(396,4524)
	(414,4550)(432,4575)(450,4599)
	(468,4624)(485,4648)(502,4671)
	(518,4694)(533,4716)(548,4737)
	(562,4757)(576,4776)(588,4794)
	(600,4812)(615,4834)(628,4856)
	(641,4877)(653,4897)(664,4917)
	(674,4936)(684,4955)(693,4973)
	(700,4991)(707,5008)(713,5025)
	(719,5040)(723,5056)(727,5070)
	(730,5084)(733,5098)(735,5111)
	(737,5124)(740,5140)(741,5156)
	(743,5172)(744,5190)(745,5211)
	(746,5233)(747,5257)(748,5283)
	(748,5311)(749,5338)(749,5363)
	(750,5384)(750,5399)(750,5408)
	(750,5411)(750,5412)
\drawline(1500,1812)(1498,1811)(1492,1809)
	(1482,1806)(1466,1800)(1445,1793)
	(1418,1783)(1385,1771)(1348,1758)
	(1307,1743)(1263,1727)(1218,1710)
	(1174,1694)(1129,1677)(1086,1660)
	(1046,1644)(1007,1629)(971,1614)
	(937,1600)(906,1586)(877,1573)
	(850,1560)(825,1548)(802,1535)
	(781,1523)(761,1511)(743,1499)
	(725,1487)(704,1471)(684,1455)
	(665,1438)(647,1421)(630,1402)
	(613,1383)(598,1363)(583,1343)
	(570,1321)(557,1299)(546,1277)
	(535,1253)(526,1230)(517,1206)
	(509,1182)(502,1158)(496,1134)
	(491,1110)(486,1085)(482,1061)
	(478,1037)(475,1012)(472,991)
	(470,969)(468,946)(466,922)
	(464,897)(463,871)(461,843)
	(460,813)(458,780)(457,746)
	(456,709)(455,671)(454,630)
	(453,589)(453,547)(452,506)
	(452,466)(451,429)(451,397)
	(450,369)(450,347)(450,331)
	(450,320)(450,315)(450,312)
\drawline(2175,2262)(2175,2261)(2175,2258)
	(2174,2249)(2173,2233)(2172,2212)
	(2170,2187)(2167,2159)(2164,2131)
	(2160,2104)(2156,2079)(2151,2056)
	(2146,2035)(2139,2015)(2131,1997)
	(2123,1979)(2113,1962)(2103,1948)
	(2093,1934)(2082,1920)(2070,1906)
	(2057,1891)(2043,1876)(2027,1861)
	(2010,1846)(1992,1831)(1974,1815)
	(1954,1800)(1934,1786)(1913,1772)
	(1892,1758)(1871,1745)(1850,1732)
	(1828,1720)(1806,1708)(1785,1697)
	(1763,1687)(1742,1678)(1721,1669)
	(1699,1660)(1677,1651)(1653,1642)
	(1630,1633)(1605,1624)(1580,1615)
	(1555,1605)(1529,1596)(1503,1587)
	(1478,1578)(1453,1569)(1428,1559)
	(1404,1550)(1381,1541)(1359,1532)
	(1338,1523)(1318,1514)(1298,1505)
	(1280,1496)(1263,1487)(1244,1477)
	(1226,1466)(1209,1454)(1192,1442)
	(1176,1429)(1161,1416)(1146,1402)
	(1132,1388)(1119,1374)(1107,1359)
	(1096,1343)(1086,1328)(1077,1313)
	(1070,1298)(1063,1283)(1058,1268)
	(1055,1254)(1052,1240)(1051,1226)
	(1050,1212)(1051,1195)(1053,1177)
	(1057,1159)(1062,1141)(1069,1122)
	(1078,1104)(1087,1084)(1098,1065)
	(1110,1047)(1122,1029)(1135,1011)
	(1148,995)(1161,979)(1174,964)
	(1187,950)(1200,937)(1213,924)
	(1226,912)(1239,900)(1253,888)
	(1267,876)(1281,864)(1295,852)
	(1309,840)(1323,829)(1336,818)
	(1348,808)(1360,798)(1371,789)
	(1382,779)(1391,771)(1400,762)
	(1411,750)(1421,739)(1431,726)
	(1441,711)(1452,695)(1463,677)
	(1474,659)(1484,641)(1492,626)
	(1497,617)(1500,613)(1500,612)
\put(3225,3762){\makebox(0,0)[lb]{{\SetFigFont{10}{12}{\rmdefault}{\mddefault}{\updefault}$\gamma_4$}}}
\put(3225,5262){\makebox(0,0)[lb]{{\SetFigFont{10}{12}{\rmdefault}{\mddefault}{\updefault}$\gamma_4^!$}}}
\put(2025,5262){\makebox(0,0)[lb]{{\SetFigFont{10}{12}{\rmdefault}{\mddefault}{\updefault}$\gamma_3^!$}}}
\put(2025,3762){\makebox(0,0)[lb]{{\SetFigFont{10}{12}{\rmdefault}{\mddefault}{\updefault}$\gamma_3$}}}
\put(2200,4587){\makebox(0,0)[lb]{{\SetFigFont{10}{12}{\rmdefault}{\mddefault}{\updefault}$p(X_\R)$}}}
\put(450,5312){\makebox(0,0)[lb]{{\SetFigFont{10}{12}{\rmdefault}{\mddefault}{\updefault}$\gamma_1^!$}}}
\put(1800,1212){\makebox(0,0)[lb]{{\SetFigFont{10}{12}{\rmdefault}{\mddefault}{\updefault}$p(X_\R)$}}}
\put(200,537){\makebox(0,0)[lb]{{\SetFigFont{10}{12}{\rmdefault}{\mddefault}{\updefault}$\gamma_1$}}}
\put(1225,1972){\makebox(0,0)[lb]{{\SetFigFont{10}{12}{\rmdefault}{\mddefault}{\updefault}$\gamma_1^!$}}}
\put(2375,2012){\makebox(0,0)[lb]{{\SetFigFont{10}{12}{\rmdefault}{\mddefault}{\updefault}$\gamma_3^!$}}}
\put(2200,412){\makebox(0,0)[lb]{{\SetFigFont{10}{12}{\rmdefault}{\mddefault}{\updefault}$\gamma_3$}}}
\put(3300,1812){\makebox(0,0)[lb]{{\SetFigFont{10}{12}{\rmdefault}{\mddefault}{\updefault}$\gamma_4^!$}}}
\put(3125,537){\makebox(0,0)[lb]{{\SetFigFont{10}{12}{\rmdefault}{\mddefault}{\updefault}$\gamma_4$}}}
\put(200,4662){\makebox(0,0)[lb]{{\SetFigFont{10}{12}{\rmdefault}{\mddefault}{\updefault}$\gamma_1$}}}
\put(850,3762){\makebox(0,0)[lb]{{\SetFigFont{10}{12}{\rmdefault}{\mddefault}{\updefault}$\gamma_2$}}}
\put(1325,4850){\makebox(0,0)[lb]{{\SetFigFont{10}{12}{\rmdefault}{\mddefault}{\updefault}$\gamma_2^!$}}}
\put(1575,325){\makebox(0,0)[lb]{{\SetFigFont{10}{12}{\rmdefault}{\mddefault}{\updefault}$\gamma_2$}}}
\put(1950,4512){\circle*{100}}
\end{picture}
}
\caption{\label{fig:real-basis}}
\end{center}
\end{figure}

Having set up the theory, it remains to produce enough examples. For this purpose we use standard approximation methods from real algebraic geometry. The outcome is as follows:

\begin{Lemma} \label{th:nash}
Let $N$ be a closed orientable manifold, equipped with a Morse function, which has the property that no two critical points lie on the same level set. One can then find: a Lefschetz fibration $p: X \rightarrow \C$ with a real structure; which comes equipped with a holomorphic volume form defined over $\R$; and a diffeomorphism $X_\R \iso N$; such that: the composition of $p_\R$ with this diffeomorphism is $C^2$-close to the given Morse function.
\end{Lemma}

\subsection{Cotangent bundles\label{subsec:cotangent}}
Fix a Lefschetz fibration with a real structure, equipped with a holomorphic volume form which is defined over $\R$. We write $N = X_\R$. By Weinstein's theorem, one can find a symplectic embedding $\kappa$ of some tubular neighbourhood of the zero-section $N \subset T^*N$ into $X$, such that $\kappa(N) = X_\R$. Because of the exactness of $X_\R$, this is an exact symplectic embedding, which means that it preserves the class of exact Lagrangian submanifolds. Since the holomorphic volume form $\eta$ is defined over $\R$, $X_\R$ admits a canonical grading, and one can use that to transfer gradings of Lagrangian submanifolds from $T^*N$ to $X$.

\begin{Lemma} \label{th:squish-1}
Let $(\tilde{L}_0,\tilde{L}_1)$ be closed admissible Lagrangian submanifolds of $T^*N$, equipped with gradings. Assume that they lie close to the zero-section, so that their images under $\kappa$ are well-defined. Then $HF^*(\kappa(\tilde{L}_0),\kappa(\tilde{L}_1)) \iso HF^*(\tilde{L}_0,\tilde{L}_1)$. Here, the Floer cohomology group on the left hand side lives in $X$, and that on the right hand side in $T^*N$.
\end{Lemma}

There is also an analogue of this for a suitable class of Lefschetz thimbles.  Namely, let $x$ be a critical point of $p$ which lies in the real locus, and $z$ its value. Let $\gamma$ be a vanishing path starting at $z$. We impose a condition similar to the one in the previous section, namely $\gamma$ should not intersect $p(X_\R)$ anywhere else, and $\im\,\gamma'(0) \neq 0$. Write $\Delta = \Delta_\gamma$ for the Lefschetz thimble, and $T^*_x \subset T^*N$ for the cotangent fibre. We fix gradings of these two submanifolds in such a way that $i(\tilde{X}_\R,\tilde\Delta;x) = 0$ and $i(\tilde{N},\tilde{T}^*_x;x) = 0$.

\begin{Lemma} \label{th:squish-2}
Let $\tilde{L}$ be a closed admissible Lagrangian submanifold of $T^*N$, with a grading. Then $HF^*(\kappa(\tilde{L}),\tilde{\Delta}) \iso HF^*(\tilde{L},\tilde{T}^*_x)$.
\end{Lemma}

\begin{Addendum} \label{th:squish-3}
Let $\tilde{L}_0,\tilde{L}_1$ be two closed admissible Lagrangian submanifolds of $T^*N$, with gradings. Then, the Floer-theoretic products in $X$ and $T^*N$ respectively,
\begin{equation}
\begin{aligned}
& HF^*(\tilde{\Delta},\kappa(\tilde{L}_1)) \otimes HF^*(\kappa(\tilde{L}_0),\tilde{\Delta}) \longrightarrow HF^*(\kappa(\tilde{L}_0),\kappa(\tilde{L}_1)), \\
& HF^*(\tilde{T}^*_x,\tilde{L}_1) \otimes HF^*(\tilde{L}_0,\tilde{T}^*_x) \longrightarrow HF^*(\tilde{L}_0,\tilde{L}_1)
\end{aligned}
\end{equation}
are compatible with the isomorphisms from Lemma \ref{th:squish-1} and \ref{th:squish-2} (together with the obvious analogue of the latter for $HF^*(\tilde\Delta,\kappa(\tilde{L}))$, which can be reduced to the original statement by duality).
\end{Addendum}

We need one more fact about gradings. Let $x$ be as before, and write
$\mu(x) = \mu(p_\R;x)$ for its Morse index as a critical point of
$p_\R$. Take two vanishing paths $\gamma$, $\gamma^!$ starting at $z =
p(x)$, with the property that $\im\,\gamma'(0) < 0$,
$\im\,(\gamma^!)'(0) > 0$. Then, the two thimbles $\Delta,\Delta^!$
and the real part $X_\R$ intersect each other pairwise  transversally at $x$.

\begin{Lemma} \label{th:triangle}
For any choice of gradings,
\begin{equation} \label{eq:minus-mu}
i(\tilde\Delta^!,\tilde\Delta;x) - i(\tilde{X}_\R,\tilde\Delta;x) - i(\tilde\Delta^!,\tilde{X}_\R;x) = -\mu(x).
\end{equation}
\end{Lemma}

\subsection{The main argument}
We start with a closed, simply-connected manifold of dimension $n \geq 6$. Choose a Morse function $q: N \rightarrow \R$ with no critical points of index $1$ or $n-1$ (it is a classical result from topology that such functions exist \cite{milnor65}). Denote by $(y_1,\dots,y_s)$ the critical points of $q$. We assume that no two such points lie on the same level set. Suppose that $(\tilde{L}_0,\tilde{L}_1)$ are closed admissible Lagrangian submanifolds of $T^*N$, equipped with gradings.
Then, {\em there is a spectral sequence converging to $HF^*(\tilde{L}_0,\tilde{L}_1)$, with}
\begin{equation} \label{eq:e1-new}
E_1^{jk} =
\Big(HF^*(\tilde{T}^*_{y_{j-r}},\tilde{L}_1) \otimes HF^*(\tilde{L}_0,\tilde{T}^*_{y_{j-r}}) \Big)^{j+k+n-\mu(y_{j-r})}.
\end{equation}
Here, we have graded the cotangent fibres as in Lemma \ref{th:squish-2}. The number  $r$ of non-real critical values appears for compatibility with labeling conventions elsewhere (it has the trivial effect of shifting the entire $E_1$ page to the right and down).

We will construct this spectral sequence by reduction to Theorem \ref{th:ss}. First, using Lemma \ref{th:nash}, find a suitable Lefschetz fibration $p: X \rightarrow \C$, whose real part approximates our Morse function; in particular, there is a bijective correspondence between the critical points of $q$ and $p_\R$, preserving Morse indices and the ordering of the critical values. As before, we extend the diffeomorphism $N \iso X_\R$ to a symplectic embedding $\kappa$ of some neighbourhood of $N \subset T^*N$ into $X$. Choose dual bases $\{\gamma_j\}$, $\{\gamma_j^!\}$ of vanishing bases which are compatible with the real structure, and take the associated Lefschetz thimbles $\{\Delta_j\}$, $\{\Delta_j^!\}$, graded according to the standard convention \eqref{eq:dual-convention}. Let's start by looking at $j \leq r$, which are just the indices corresponding to non-real critical values. In that case, $\Delta_j$ is disjoint from $X_\R$, hence also from $\kappa(L_0)$, $\kappa(L_1)$ if we bring those submanifolds sufficiently close to the zero-section. As a consequence, the associated columns in \eqref{eq:e1} vanish. Now consider the real critical points $x_j$ and critical values $z_j$, $j > r$. By Lemma \ref{th:triangle}, there is some constant $d_j$ such that
\begin{equation}
\begin{aligned}
 & i(\tilde{X}_\R,\tilde\Delta_j;x_j) = d_j, \\
 & i(\tilde{X}_\R,\tilde\Delta_j^!;x_j) = n - i(\tilde\Delta_j^!,\tilde{X}_\R;x_j) = n - \mu(x_j) + d_j.
\end{aligned}
\end{equation}
Lemma \ref{th:squish-2}, as originally formulated, assumes that these two quantities are zero. We adjust the statement to take into account the gradings here, and find that
\begin{equation}
\begin{aligned}
 & HF^*(\kappa(\tilde{L}_0),\tilde\Delta_j) \iso HF^{*-d_j}(\tilde{L}_0,\tilde{T}^*_{x_j}), \\
 & HF^*(\tilde\Delta_j^!,\kappa(\tilde{L}_1))
 \iso HF^{*+n-\mu(x_j)+d_j}(\tilde{T}^*_{x_j},\tilde{L}_1).
\end{aligned}
\end{equation}
In view of the correspondence between critical points of $p_\R$ and $q$, and the isotopy invariance of Floer cohomology in $T^*N$ (which makes it irrelevant whether one takes the cotangent fibre at $x_j$ or $y_{j-r}$), the $E_1$ page now takes on the form \eqref{eq:e1-new}. In fact, the last-mentioned observation also shows that all columns on this page are isomorphic, up to a shift.

For our first application, take $\tilde{L}_0 = \tilde{L}_1 = \tilde{L}$ equal, and write
\begin{equation}
H = HF^*(\tilde{T}^*_x,\tilde{L}) \otimes HF^*(\tilde{L},\tilde{T}^*_x),
\end{equation}
where $x$ is any point in $N$. Let $a \leq b$ be the lowest and highest degrees in which $H$ is nonzero. Then, the term in \eqref{eq:e1-new} with the highest total (row plus column) degree is $E_1^{m,b-m}$ (recall that by convention, the last critical value is the maximum of $p_\R$, hence the unique critical point with Morse index $n$). Moreover, because there are no critical points with index $n-1$, any other term has total degree $\leq b-2$, so this highest degree piece necessarily survives to $E_\infty$. The lowest degree piece, which is $E_1^{r+1,a-n-r-1}$, survives for analogous reasons. Now, the spectral sequence converges to $HF^*(\tilde{L},\tilde{L}) = H^*(L;\K)$, which is obviously concentrated in degrees $0 \leq * \leq n$. Hence, $b \leq n$ and $a-n \geq 0$, which implies equality. Moreover, the bottom and top classical cohomology groups are one-dimensional, hence we find that {\em $HF^*(\tilde{L},\tilde{T}^*_x)$ is one-dimensional}. After shifting the grading of $L$, we may assume that it is concentrated in degree $0$. Of course, by passing to Euler characteristics, it follows that projection $L \rightarrow N$ has degree $\pm 1$.

Next, apply the same spectral sequence to $\tilde{L}_0 = \tilde{L}$
and $\tilde{L}_1 = \tilde{N}$ (at this point, we have to assume that $N$ is spin). The top degree piece in the $E_1$ page is $E_1^{m,n-m} = HF^n(\tilde{T}^*_x,\tilde{N}) \otimes HF^0(\tilde{L},\tilde{T}^*_x) \iso \K$, and this survives for the same reasons as before. As a consequence, the edge homomorphism is an isomorphism in degree $n$. Addendum \ref{th:edge} describes the edge homomorphism as a Floer-theoretic product in $X$, and we can apply Addendum \ref{th:squish-3} to transfer the product to $T^*N$, where it takes the form $HF^n(\tilde{T}^*_x,\tilde{N}) \otimes HF^0(\tilde{L},\tilde{T}^*_x) \rightarrow HF^n(\tilde{L},\tilde{N})$. Here, $x$ can be arbitrary by isotopy invariance. Transferring this insight back to $X$, we find that {\em for any $j>r$, the map
\begin{equation}
 HF^n(\tilde\Delta_j,\kappa(\tilde N)) \otimes HF^0(\kappa(\tilde{L}),\tilde\Delta_j)
 \longrightarrow HF^n(\kappa(\tilde{L}),\kappa(\tilde{N}))
\end{equation}
is an isomorphism of one-dimensional vector spaces}. Equivalently, fixing a nonzero element $c$ in the dual group $HF^0(\kappa(\tilde{N}),\kappa(\tilde{L}))$, the statement is that product with $c$ induces isomorphisms $HF^*(\kappa(\tilde L),\tilde\Delta_j) \rightarrow HF^*(\kappa(\tilde N),\tilde\Delta_j)$ for all $j>r$. In this form, the statement also holds for $j \leq r$, where the groups involved are all zero. Hence, Addendum \ref{th:cone} applies. Recall that (closed) admissible Lagrangians and their Floer cohomology groups form a (genuine) category $H\F(X)$.  We may therefore appeal to the Yoneda Lemma, which asserts that in any category $C$, an object $O$ is determined up to isomorphism by the functor $\star \mapsto \textrm{Mor}_{C}(O, \star)$. The conclusion of the Addendum then implies that $c$ is an isomorphism in the category $H\F(X)$. In particular, $\kappa(\tilde{L})$ and $\kappa(\tilde{N})$ must have isomorphic endomorphism rings, so $H^*(L;\K) \iso HF^*(\kappa(\tilde{L}),\kappa(\tilde{L})) \iso HF^*(\kappa(\tilde{N}),\kappa(\tilde{N})) \iso H^*(N;\K)$. We already know that the projection $L \rightarrow N$ has degree $\pm 1$, hence is injective on cohomology. By comparing dimensions, it follows that it must be an isomorphism. Now let $L_0,L_1$ be two submanifolds as in the last part of Theorem \ref{th:main}. We know that in $H\F(X)$, both $\kappa(\tilde{L}_j)$ are isomorphic to $\kappa(\tilde{N})$, hence $HF^*(\kappa(\tilde{L}_0),\kappa(\tilde{L}_1)) \iso HF^*(\kappa(\tilde{N}),\kappa(\tilde{N})) = H^*(N;\K)$. By definition of the Floer cochain complex, this implies the desired lower bound on the number of intersection points. Finally, to get rid of the assumption $\mathrm{dim}(N) \geq 6$, one argues by stabilization (taking the product with a sphere of large dimension).

\section{Parallel transport}

This section contains the proof of Lemma \ref{th:parallel}. The point is to show that the parallel transport vector fields are integrable, which means that orbits do not escape to infinity in finite time.

\subsection{Estimates for horizontal vector fields}
Take the vector field $\partial_z$ on the base $\C$, and take its unique lift to a horizontal vector field on $X$ (away from the critical points), using the K{\"a}hler metric. The outcome can be written as
\begin{equation} \label{eq:horiz}
\xi = \frac{\nabla p}{||\nabla p||^2}.
\end{equation}
Let's place ourselves at a point at infinity in $\bar{X}$, which lies in the closure of $s^{-1}(0)$ (note that this is equal to the closure of $p^{-1}(z)$, for any $z$). There are local coordinates $(x_1,\dots,x_n)$ around that point, and a local trivialization of $\O_{\bar{X}}(1)$, with respect to which
\begin{equation}
s(x) = x_{k+1}, \quad t(x) = x_1^{m_1} \cdots x_k^{m_k}.
\end{equation}
We then have
%
\begin{align} \label{eq:1-bound}
 & ||\nabla p|| \gapprox |\partial_j p| \gapprox
 \frac{|p|}{|x_j|} \text{ for $1 \leq j \leq k$; and}\\
 & ||\nabla t|| \lapprox \sum_{j=1}^k \frac{|t|}{|x_j|}. \label{eq:2-bound}
\end{align}
Here, the gradient and its norm are formed with respect to the given K{\"a}hler metric. The notation $\gapprox$ means that inequality holds up to some multiplicative constant (mainly, this involves comparing the metric with the standard one). In our local trivialization, the hermitian metric on $\O_{\bar{X}}(1)$ is $||\cdot||^2 = e^\sigma |\cdot|^2$ for some smooth function $\sigma$, hence $h = -\log ||t||^2 = -\log |t|^2 - \sigma$. Differentiate this in direction of \eqref{eq:horiz}:
\begin{equation}
\begin{aligned}
 \big| \xi.h \big| & \leq \frac{|\leftsc \nabla p,\nabla \sigma\rightsc |}{||\nabla p||^2} +
 \frac{2 |t| \cdot |\leftsc \nabla t, \nabla p \rightsc|}{||\nabla p||^2 \cdot |t|^2} \\ &
 \lapprox \frac{1}{||\nabla p||} + \frac{||\nabla t||}{||\nabla p|| \cdot |t|}.
\end{aligned}
\end{equation}
In view of \eqref{eq:1-bound}, the first term is bounded above by $const/|p|$ near $x = 0$. A combination of \eqref{eq:1-bound} and \eqref{eq:2-bound} yields the same bound on the second term.

\subsection{Application}
Consider parallel transport along a horizontal segment $[a,b] \subset \C$, which avoids all critical values. This is defined by integrating $\xi$ over $p^{-1}([a,b])$. Assume temporarily that $0 \notin [a,b]$, so that we get a bound on $1/|p|$. Then, after covering the closure of $s^{-1}(0)$ with finitely many neighbourhoods of the kind considered above, it follows that $|\xi.h|$ is bounded on the whole of $p^{-1}([a,b])$. Since $h$ is an exhausting function, that gives an a priori bound on the growth of trajectories, prohibiting their escape to infinity.
To get rid of the assumption $0 \notin [a,b]$, one argues as follows: if $p$ is a Lefschetz fibration, then so is $p+c$ for any constant $c$. Moreover, the parallel transport maps remain the same. In other words, by changing the way in which we choose local coordinates, the $1/|p|$ bound can be replaced by a $1/|p+c|$ one. It then suffices to choose $c \notin [a,b]$.

Finally, parallel transport along an arbitrary path $\beta$ is defined by taking the horizontal lifts of the tangent vectors $\partial\beta/\partial r$. Since these are all complex multiples of \eqref{eq:horiz}, the same argument as before works.

\section{From Morse functions to Lefschetz fibrations}

This section contains the proof of Lemma \ref{th:nash}. This is a standard exercise in real algebraic geometry, using the Nash-Tognioli approximation theorem, resolution of singularities, and Bertini-type transversality results.

\subsection{Real algebraic approximation}
Let $N$ be a closed $n$-dimensional manifold, smoothly embedded into $\R^{2n+1}$. The Nash-Tognioli theorem, in the form given in \cite[Theorem 7]{ivanov82}, says that any such $N$ can be $C^\nu$-approximated by a smooth real algebraic variety $U_\R$, for every $\nu$ ($\nu = 2$ will be enough for us). Let $U$ be the complexification of $U_\R$, which means the affine variety in $\C^{2n+1}$ defined by the same equations as $U_\R$. In general, this complexification will be singular, but one can throw out the singularities by using the following trick \cite{kulkarni78}. $Sing(U) \subset U$ is itself an algebraic subvariety defined over $\R$, which means the zero-set of real polynomials $f_1,\dots,f_r$. Set $F = f_1^2 + \cdots + f_r^2$, and consider
\begin{equation} \label{eq:graph-1}
\{(u,v) \in \C^{2n+1} \times \C \;:\; u \in U, \; F(u)v = 1\}.
\end{equation}
This is isomorphic to $U \setminus Sing(U)$. Since $U_\R$ is smooth, we have $Sing(U) \cap U_\R = \emptyset$, which means that the real part of \eqref{eq:graph-1} is isomorphic to $U_\R$ (in both cases, the isomorphism is given by projection to $\C^{2n+1}$).

From now on, we assume that $N$ is orientable. Take an $n$-form $\beta_\R$ on $\R^{2n+1}$ (not closed, of course) whose restriction to $N$ is a volume form. We may choose $\beta_\R$ to be real algebraic; this is just the Stone-Weierstrass theorem on polynomial approximation. In our previous construction, we take $U_\R$ to be sufficiently close to $N$, and then $\beta_\R|U_\R$ will again be a volume form. Take the complexification $\beta$, and pull it back to \eqref{eq:graph-1} by projection. This may not necessarily be a complex volume form, but the set where it degenerates (becomes zero) is an algebraic subvariety defined over $\R$, and disjoint from the real locus. We choose real defining polynomials $g_1,\dots,g_s$ for this subvariety, and apply the same trick as before, which means passing to
\begin{equation} \label{eq:graph-2}
\{(u,v,w) \in \C^{2n+1} \times \C \times \C \; : \; u \in U,\; F(u)v = 1, \; G(u,v)w = 1\},
\end{equation}
where $G = g_1^2 + \cdots + g_s^2$. The outcome is that we have a smooth affine algebraic variety defined over $\R$, whose real part is diffeomorphic to $N$, which comes equipped with a holomorphic volume form.

\subsection{Resolution of singularities at infinity}
The projective closure of \eqref{eq:graph-2} is not in general well-behaved (it can have arbitrarily bad singularities). To resolve these, we appeal to Hironaka's theorem \cite{hironaka64}. The precise statement is as follows: for some $q \gg 2n+3$, there is an affine algebraic variety $X \subset \C^q$, such that:
\begin{itemize} \itemsep0.5em
\item $X$ is defined over $\R$;
\item projection to the first $2n+3$ coordinates maps $X$ isomorphically to \eqref{eq:graph-2};
\item the projective closure $\bar{X} \subset \CP{q}$ is smooth, and the divisor at infinity $Y = X \cap \CP{q-1}$ has at most normal crossing singularities.
\end{itemize}
We equip $\bar{X}$ with the line bundle $\O_{\bar{X}}(1)$, with its standard metric, and the section $t$ defining the divisor at infinity. Finally, we take the previously constructed holomorphic volume form, and pull it back to a form $\eta$ on $X$, then extend that to a rational form on $\bar{X}$. This data satisfies all the conditions from Sections \ref{subsec:lefschetz} and \ref{subsec:squish}, as far as the geometry of the total space $X$ itself is concerned. Next, we will address the construction of the real Lefschetz pencil $p$.

\subsection{Constructing the Lefschetz fibration}
Suppose that our $N$ comes with a choice of Morse function, which has at most one critical point in each level set. At the outset of the construction, we may assume that the embedding $N \subset \R^{2n+1}$ has been chosen in such a way that the first coordinate $u_1$ is $C^2$-close to the given Morse function (one can even arrange that the two are equal, but we won't need this). By choosing $U_\R$ sufficiently close to $N$ and going through the construction, one gets the following: there is a section $s$ of $\O_{\bar{X}}(1)$, which is defined over $\R$, such that the restriction of $p = s/t$ to the real part $X_\R \iso N$ is $C^2$-close to the original Morse function.

By Bertini's theorem, the complex hyperplanes which intersect $\bar{X}$ non-transversally form a proper subvariety of the dual projective space $\CP{q-1}$. The real locus of that is a proper subvariety of $\RP{q-1}$, hence its complement is open and dense. This means that by a small perturbation of $s$ inside the space of sections defined over $\R$, we may achieve that $s^{-1}(0)$ is smooth. Similarly, a generic choice ensures that $s$ intersects all the strata of $Y$ transversally; that the critical points of $p$ are nondegenerate; and that at most one such point lies in each fibre. In all those cases, it is a classical fact that the set of complex parameter values (choices of $s$) where things go wrong is a constructible subset of positive codimension, and one applies the same argument as before to obtain the desired result for real $s$.

\section{Grading issues}

In this section, we review in a little more detail the standard machinery of graded Lagrangian submanifolds, and the resulting gradings on Floer cohomology groups. Most of our argument, including Lemma \ref{lem:gradingconvention}, uses this machinery only in straightforward ways. The exception is Lemma \ref{th:triangle}, for which we will provide a proof based on the index formula for holomorphic triangles (more pedestrian proofs, by explicit computation of all the indices involved, are also possible).

\subsection{Generalities}
Let $M$ be any symplectic manifold, and $Gr = Gr_M \longrightarrow M$ the bundle of Lagrangian Grassmannians associated to the symplectic vector bundle $TM$. Suppose that in addition, we have an infinite cyclic covering
\begin{equation} \label{eq:tilde}
\widetilde{Gr} \longrightarrow Gr,
\end{equation}
which fibrewise is isomorphic to the universal covering of each
Grassmannian. Every Lagrangian submanifold $L \subset M$ comes with a
tautological section of $Gr|L$, given by $x \mapsto TL_x$. One defines
a grading of $L$ to be a lift of this to $\widetilde{Gr}$, which means
a choice of preimage  $\widetilde{TL}_x \in \widetilde{Gr}_x$ for any $x$, varying continuously. A graded Lagrangian submanifold $\tilde{L}$ is a Lagrangian submanifold equipped with a choice of grading. It is then obvious that the Maslov class $m_L \in H^1(L)$, defined as the pullback of the element of $H^1(Gr)$ classifying \eqref{eq:tilde} by the tautological section, is the obstruction to the existence of a grading.

One case where this formalism applies in a straightforward way is that of cotangent bundles $T^*N$. Take the tautological section associated to the zero-section $N \subset T^*N$, and extend that over the whole of $T^*N$ in an arbitrary way. Then, take $\widetilde{Gr}_{T^*N}$ to be the fibrewise universal cover with base points given by that section. As a direct consequence of the definition, the zero-section $N$ comes with a trivial grading (if one thinks of points in the universal cover as equivalence classes of paths, this is given for each $x$ by the constant path at $TN_x$). In a more general context, this example occurs as follows: suppose that $M$ is any symplectic manifold equipped with a covering \eqref{eq:tilde}, and $\tilde{N}$ a graded Lagrangian submanifold. Enlarge the inclusion $N \hookrightarrow M$ to a symplectic embedding of a tubular neighbourhood of the zero-section inside $T^*N$. On this neighbourhood, there is a preferred isomorphism between the pullback of $\widetilde{Gr}_M$ and the previously considered covering $\widetilde{Gr}_{T^*N}$, given by the grading of $N \subset M$. In less precise but more practical terminology, there is a unique coherent way of mapping graded Lagrangian submanifolds in $T^*N$ to ones in $M$, with the property that the zero-section (with its trivial grading) gets mapped to $\tilde{N}$.

Algebro-geometrically, the natural source of coverings \eqref{eq:tilde} is as follows. Let $X$ be a K{\"a}hler manifold with a holomorphic volume form $\eta$. This induces a squared phase function $\alpha: Gr_X \rightarrow S^1$, defined by $\alpha(\Lambda) = \eta(v_1 \wedge \cdots \wedge v_n)^2/|\eta(v_1 \wedge \cdots \wedge v_n)|^2$, where $\{v_j\}$ is any basis of the Lagrangian subspace $\Lambda \subset TX_x$. Taking \eqref{eq:tilde} to be the pullback of $\R \rightarrow S^1$ by $\alpha$, one immediately sees that: for any Lagrangian submanifold $L$, the class $m_L$ is represented by the function $\alpha_L: L \rightarrow S^1$, $\alpha_L(x) = \alpha(TL_x)$; and a grading of $L$ is the same as a real-valued phase function $\tilde\alpha_L: L \rightarrow \R$ satisfying $\exp(2\pi i \tilde\alpha_L) = \alpha_L$. In our specific application, we have a real involution (which reverses the K{\"a}hler form and almost complex structure, and maps $\eta$ to $\bar\eta$). In that case, $\eta(v_1 \wedge \cdots \wedge v_n) \in \R$ for any basis $\{v_j\}$ of $(TX_\R)_x$, hence $\alpha_{X_\R} \equiv 1$. One therefore has a canonical grading of the real locus, $\tilde\alpha_{X_\R} \equiv 0$.

\subsection{Maslov indices}
The main role of gradings is to allow us to fix the $\Z$-grading of the Floer cochain complex. This is done through the absolute Maslov index, which was already mentioned in \eqref{eq:absolute} above. We refer to \cite{seidel99} for a general definition. For our applications, only one case is really important. Take $T^*N$ with its standard covering $\widetilde{Gr}$. Suppose that $L_0 = N$ is the zero-section, and $L_1 = \mathrm{graph}(df)$ the graph of some exact one-form $df$. There is an obvious isotopy from $L_0$ to $L_1$, and we assume that gradings have been chosen in a way that is compatible with this isotopy. Then, for any nondegenerate critical point of $f$, the absolute Maslov index equals the Morse index:
\begin{equation} \label{eq:i-equal-m}
i(x) = \mu(x).
\end{equation}

We now prove Lemma \ref{lem:gradingconvention}, which holds for any Lefschetz fibration, namely that \eqref{eq:dual-convention} is compatible with the isotopy $\Delta_m \htp \Delta_m^!$.  Pick local co-ordinates near the critical point for which $p(x) = x_1^2+\cdots+x_n^2$.
We first move the relevant vanishing paths so that $\gamma_m(r) = z_m + \exp(-i\epsilon)r$ and $\gamma_m^!(r) = z_m + \exp(i\epsilon)r$ near $r = 0$, for some small $\epsilon>0$. It is well-known, see for instance \cite[Lemma 1.7]{seidel01}, that one can deform the K{\"a}hler structure locally to make it standard in any given complex coordinate system, so there are local coordinates in which the vanishing cycles are $\Delta_m = \exp(-i\epsilon/2)\R^n$, $\Delta^!_m = \exp(i\epsilon/2)\R^n$. Rotate linearly into Darboux coordinates $(p,q)$, with symplectic form $dp\wedge dq$, so that $\Delta^!_m = \{p = 0\}$ and $\Delta_m = \{p = \tan(\epsilon)q\}$. The isotopy $\Delta_m \htp \Delta_m^!$ is locally given by deforming $\epsilon$ to zero, and one can apply \eqref{eq:i-equal-m} to show that if one chooses gradings compatibly with this isotopy, the Maslov index is indeed zero.  By contrast, in the analogous local model describing $\Delta_1 = \{p=-\tan(\epsilon)q\}$ as a graph over $\Delta^!_1=\{p=0\}$, the generating function has a maximum (rather than a minimum), which gives a discrepancy between \eqref{eq:i-equal-m} and \eqref{eq:dual-convention}.

\subsection{An index computation}
Next, let $L_0,L_1,L_2$ be three graded Lagrangian submanifolds, intersecting at the same point $x$, with pairwise transverse intersections. Once one has chosen a compatible almost complex structure, there is a trivial holomorphic triangle $u$ with boundary conditions (anticlockwise ordered) $L_0,L_1,L_2$, namely the constant $u \equiv x$. Let $D_u$ be the linearized operator at $u$. Equivalently, this is the $\bar\partial$-operator on the trivial vector bundle with fibre $TM_x$ over a three-punctured disc, with boundary values in $TL_{0,x}$, $TL_{1,x}$, $TL_{2,x}$. As a special case of the general index formula for such operators (one reference with compatible terminology is \cite[Proposition 11.13]{seidel04}, but there are many others), one has the following equality: for any choice of gradings,
\begin{equation} \label{eq:triangle}
\index(D_u) = i(\tilde{L}_0,\tilde{L}_2;x) - i(\tilde{L}_0,\tilde{L}_1;x) - i(\tilde{L}_1,\tilde{L}_2;x).
\end{equation}
The simplest case is that of three lines in the plane. Thinking of this as $\C = T^*\R$, we can compute absolute Maslov indices (for suitable choices of gradings) from \eqref{eq:i-equal-m}, and then apply \eqref{eq:triangle}. The outcome is that the index of $D_u$ is either $0$ or $-1$, depending on the ordering of the lines (Figure \ref{fig:2-triangles}).

\begin{figure}
\setlength{\unitlength}{3947sp}%
\begingroup\makeatletter\ifx\SetFigFont\undefined%
\gdef\SetFigFont#1#2#3#4#5{%
  \reset@font\fontsize{#1}{#2pt}%
  \fontfamily{#3}\fontseries{#4}\fontshape{#5}%
  \selectfont}%
\fi\endgroup%
\begin{picture}(3837,1861)(226,-1310)
\put(2520,-1261){\makebox(0,0)[lb]{\smash{{\SetFigFont{10}{12.0}{\rmdefault}{\mddefault}{\updefault}{Case 2: $\index D_u = 0$}%
}}}}
\thinlines
{\put(601,-811){\line( 2, 3){900}}
}%
{\put(301,-136){\line( 1, 0){1500}}
}%
{\put(2851,539){\line( 2,-3){900}}
}%
{\put(2851,-811){\line( 2, 3){900}}
}%
{\put(2551,-136){\line( 1, 0){1500}}
}%
\put(226,-306){\makebox(0,0)[lb]{\smash{{\SetFigFont{10}{12.0}{\rmdefault}{\mddefault}{\updefault}{$L_0$}%
}}}}
\put(451,-961){\makebox(0,0)[lb]{\smash{{\SetFigFont{10}{12.0}{\rmdefault}{\mddefault}{\updefault}{$L_1$}%
}}}}
\put(1351,-961){\makebox(0,0)[lb]{\smash{{\SetFigFont{10}{12.0}{\rmdefault}{\mddefault}{\updefault}{$L_2$}%
}}}}
\put(2476,-306){\makebox(0,0)[lb]{\smash{{\SetFigFont{10}{12.0}{\rmdefault}{\mddefault}{\updefault}{$L_0$}%
}}}}
\put(2701,-961){\makebox(0,0)[lb]{\smash{{\SetFigFont{10}{12.0}{\rmdefault}{\mddefault}{\updefault}{$L_2$}%
}}}}
\put(3601,-961){\makebox(0,0)[lb]{\smash{{\SetFigFont{10}{12.0}{\rmdefault}{\mddefault}{\updefault}{$L_1$}%
}}}}
\put(246,-1261){\makebox(0,0)[lb]{\smash{{\SetFigFont{10}{12.0}{\rmdefault}{\mddefault}{\updefault}{Case 1: $\index D_u = -1$}%
}}}}
{\put(601,539){\line( 2,-3){900}}
}%
\end{picture}%
\caption{\label{fig:2-triangles}}
\end{figure}

We now turn to Lemma \ref{th:triangle}. First consider a simplified local model, namely $X = \C^n$ with the standard symplectic structure, the constant complex volume form $\eta = dx_1 \wedge \cdots \wedge dx_n$, and the standard real structure (complex conjugation).  Moreover, our $p$ should be a quadratic function
\begin{equation} \label{eq:p-standard}
p(x) = -x_1^2 - \cdots - x_\mu^2 + x_{\mu+1}^2 + \cdots + x_n^2,
\end{equation}
and our vanishing paths are straight lines $\gamma(r) = -ir$, $\gamma^!(r) = ir$. The associated Lefschetz thimbles are $\Delta = \sqrt{i}\R^\mu \times \sqrt{-i}\R^{n-\mu}$, $\Delta^! = \sqrt{-i}\R^\mu \times \sqrt{i}\R^{n-\mu}$. From \eqref{eq:triangle} we know that the right hand side of \eqref{eq:minus-mu} is the index of the linearized operator $D_u$ for the constant holomorphic triangle $u \equiv 0$ with boundary conditions $(\Delta^!,X_\R,\Delta)$. Clearly, this operator splits into the direct sum of $n$ scalar ones, of which the first $\mu$ correspond to the left-hand picture in Figure \ref{fig:2-triangles}, and the remaining ones to the right-hand one. Hence, its index is $-\mu$ by our previous computation.

To derive the general case from this, note first that by the real-analytic version of the Morse Lemma, there are always coordinates compatible with the real structure in which, near a critical point $x$ of $p_{\R}$ of Morse index $\mu$, $p$ has the form \eqref{eq:p-standard}.  We apply a local deformation of the K\"ahler form as before;
in our case, this can be done compatibly with the real structure, so that $X_\R$ remains Lagrangian. Obviously, such a deformation also affects $\Delta$ and $\Delta^!$, which change by a Lagrangian isotopy. However, throughout this isotopy they do stay transverse to each other, as well as to $X_\R$, so the relevant indices remain the same.

\section{A shrinking argument}

We will now prove the remaining results from Section \ref{subsec:cotangent} (everything other than Lemma \ref{th:triangle}). The idea is to arrange that the relevant holomorphic curves have very small energy, and therefore cannot escape a neighbourhood of the real locus $X_\R$.

\subsection{The Monotonicity Lemma}
Let $(M,\o,J)$ be any compact symplectic manifold, equipped with a compatible almost complex structure, and $L \subset M$ a Lagrangian submanifold. The Monotonicity Lemma, in its relative form, says:

\begin{Lemma} \label{th:monotonicity}
There are constants $\rho > 0$, $\gamma > 0$ such that the following holds. Let $B = B(r;y) \subset M$ be a closed ball (in the associated Riemannian metric) of radius $0<r \leq \rho$ around a point $y$. Let $\Sigma$ be a compact connected Riemann surface with corners, and $u: \Sigma \rightarrow B$ a non-constant $J$-holomorphic curve such that $y \in u(\Sigma)$, $u(\partial \Sigma) \subset L \cup \partial B$. Then
\begin{equation} \label{eq:m-energy}
E(u) = \int_\Sigma u^*\o \geq \gamma r^2.
\end{equation}
\end{Lemma}

The case where $L = \emptyset$ is the most familiar one, see for instance \cite{audin-lafontaine}, but the relative version is proved in the same way (by looking at Darboux charts and doing integration by parts).
Now take a closed manifold $N$, fix a Riemannian metric on it, and denote by $V = T^*_{< \lambda}N$ the subspace of cotangent vectors of length less than some constant $\lambda$ (similarly, we write $T^*_{>\lambda}N$ for the set of vectors of length $>\lambda$). This comes with a canonical symplectic form $\o_V = d\theta_V$ and (due to the metric) almost complex structure $J_V$. Fix some point $x \in N$, and let $\Delta_V = T^*_{<\lambda,x} \subset V$ be the piece of its cotangent fibre which is contained in $V$.

\begin{Lemma} \label{th:pass}
For fixed $\lambda$, there is a constant $\epsilon>0$ with the following property. Let $\Sigma$ be a compact connected Riemann surface with corners, and $u: \Sigma \rightarrow V$ a $J_V$-holomorphic map with
\begin{equation}
u(\partial\Sigma) \subset \Delta_V \cup T^*_{<\lambda/3} \cup T^*_{>2\lambda/3}.
\end{equation}
Suppose that the image $u(\Sigma)$ contains some cotangent vector of length $< \lambda/3$, and another one of length $> 2\lambda/3$. Then $E(u) \geq \epsilon$.
\end{Lemma}

This can be derived directly from the previous result, as follows. By connectedness, we know that $u(\Sigma)$ must contain a $y$ with $||y|| = \lambda/2$. There is some $r$ depending only on the metric, such that the ball $B(r;y)$ is disjoint from $T^*_{<\lambda/3}$ and $T^*_{>2\lambda/3}$. By making $r$ smaller if necessary, we may assume that it is less than the constant $\rho$ from Lemma \ref{th:monotonicity}. We may also find some $r/2 \leq r' \leq r$, such that $u$ is transverse to the boundary of $B(r';y)$. Setting $\Sigma' = u^{-1}(B(r';y))$ and applying \eqref{eq:m-energy} to $u' = u|\Sigma'$, we find that $E(u) \geq E(u') \geq \gamma (r')^2 \geq \gamma r^2/4 = \epsilon$.

\subsection{Energy estimates}
We now turn to our application. Supposing that $\lambda$ is sufficiently small, we have a symplectic embedding $\kappa: V \longrightarrow X$ such that $\kappa|N$ is the given identification of $N$ with $X_\R$. Moreover, this embedding is exact, which means that $\kappa^*\theta$ differs from the canonical one-form $\theta_V$ by an exact one-form (this is obvious, because both $\kappa^*\theta$ and $\theta_V$ vanish on the zero-section). After making a suitable change
\begin{equation} \label{eq:change-theta}
\theta \longmapsto \theta + dH,
\end{equation}
and maybe shrinking $\lambda$ a little, we may assume that $\kappa^*\theta = \theta_V$. Such a change is unproblematic, because our only use for $\theta$ is in defining the class of exact Lagrangian submanifolds, which is unaffected by \eqref{eq:change-theta}. Similarly, we can find a compatible almost complex structure $J$ on $X$ which agrees with the given complex structure outside a compact subset, and such that $\kappa^*J = J_V$.

Let $(L_0,L_1)$ be two Lagrangian submanifolds of $V$ which are closed and exact, and which intersect transversally. Exactness means that there are functions $K_0,K_1$ such that $dK_j = \theta_V|L_j$. Recall that the action functional at an intersection point $x \in L_0 \cap L_1$ is defined to be
\begin{equation} \label{eq:action}
A(x) = A_{L_0,L_1}(x) = K_1(x) - K_0(x).
\end{equation}
In our case, this agrees with the action functional for $\kappa(L_0)$, $\kappa(L_1)$ viewed as Lagrangian submanifolds of $X$. In particular, if $u: \R \times [0;1] \rightarrow X$ is any $J$-holomorphic strip with $u(\R \times \{j\}) \subset \kappa(L_j)$ and limits $\lim_{s \rightarrow \pm \infty} u(s,\cdot) = x_{\pm}$, then
\begin{equation} \label{eq:energy}
E(u) = A(x_-) - A(x_+).
\end{equation}
Now rescale our Lagrangian submanifolds radially in the cotangent bundle, replacing them by $\rho L_j$ for some $0<\rho<1$. Because $\theta_V$ is homogeneous in fibre direction, the associated functions change to $\rho K_j$, which by \eqref{eq:action} means that the action is $A_{\rho L_0,\rho L_1}(\rho x) = \rho A(x)$. Using \eqref{eq:energy} we see that there is a constant $C$ such that $E(u) \leq C\rho$ for any finite energy holomorphic strip $u$ with boundary in $(\kappa(\rho L_0),\kappa(\rho L_1))$. By choosing $\rho$ sufficiently small, one can arrange that $(\rho L_0,\rho L_1)$ are contained in $T^*_{<\lambda/3} N$, and that $C\rho$ is less than the constant $\epsilon$ from Lemma \eqref{th:pass} (for this particular application, the cotangent fibre $\Delta_V$ is irrelevant). It follows that no such $u$ can leave $T^*_{<2\lambda/3}N$. On the other hand, the Floer cohomology of $(\rho L_0,\rho L_1)$ in $T^*_{<2\lambda/3}N$ is the same as that in the entire cotangent bundle, by the maximum principle. This essentially completes the proof of Lemma \ref{th:squish-1}. A little caution must be observed, since Floer cohomology is usually computed by using pseudo-holomorphic strips for a generic $t$-dependent perturbation of the almost complex structure. However, one can make this perturbation small and supported inside $T^*_{<\lambda/3}N$, and then the pseudo-holomorphic strips must still remain inside $T^*_{<2\lambda/3}N$, by Gromov compactness.

\subsection{The Lefschetz thimble case}
Let $\Delta \subset X$ be a Lefschetz thimble for a path chosen as in Lemma \ref{th:squish-2}. In particular, $\Delta$ intersects $X_\R$ in a single point $x$, and the intersection is transverse there. Let $\Delta_V \subset V$ be the cotangent fibre at the same point. One can find local Darboux coordinates $(p,q)$ for $X$ centered at $x$, in which
\begin{equation}
 X_\R = \{q = 0\}, \quad \kappa(\Delta_V) = \{p = 0\}, \quad \Delta = \{p = df(q)\}
\end{equation}
for some function $f$ which has a critical point at $q = 0$. By multiplying $f$ with a smooth cutoff function vanishing near $q = 0$, one finds another Lagrangian submanifold $\Delta' \subset X$ isotopic to $\Delta$, such that $\Delta' \cap X_\R = \{x\}$ and $\Delta' = \kappa(\Delta_V)$ near $x$. After making $\lambda$ smaller if necessary, we may assume that in fact, $\kappa^{-1}(\Delta') = \Delta_V$.

Now apply the previous argument to $(\kappa(L_0),\Delta')$, where $L_0 \subset V$ is closed and exact.
Since $\theta_V|\Delta_V = 0$, we may choose a function $K$ with $dK = \theta|\Delta'$ in such a way that it vanishes on $\kappa(\Delta_V)$. In that case, a linear rescaling of $L_0$ still results in a linear change of the actions of all intersection points, hence in a linear decrease of the energy. One applies Lemma \ref{th:pass} as before, to argue that all holomorphic strips must remain inside $V$, and obtains Lemma \ref{th:squish-2}. Addendum \ref{th:squish-3} is the same argument applied to holomorphic triangles; this works because the energy of a holomorphic triangle is determined by the actions of its endpoints, in a way which is entirely parallel to \eqref{eq:energy}.

\section{The spectral sequence}

This section concerns Theorem \ref{th:ss}. This is essentially the same as \cite[Corollary 18.27]{seidel04}, but a review still seems appropriate, if only because we need to derive some additional properties of the spectral sequence. The idea originally arose in algebraic geometry, where the prototype is Beilinson's spectral sequence for sheaves on projective space. This was subsequently generalized to triangulated categories admitting full exceptional collections, see \cite[Section 2.7.3]{gorodentsev-kuleshov04}. We follow this partially, but combine it with a more direct approach in terms of $A_\infty$-modules. Having derived the spectral sequence in this purely algebraic framework, we then quote, without proof, the geometric results from \cite[Chapter 3]{seidel04} which explain how this applies to Lefschetz fibrations.  An informal overview of the geometric side of the story is also given in \cite{fukaya-seidel-smith07b}.


\subsection{$A_\infty$-modules}

Let $\A$ be a directed $A_\infty$-category, linear over $\K$. By definition (see \cite{seidel-mu} or \cite[Section 5m]{seidel04}) this is a strictly unital $A_\infty$-category, with a finite ordered set of objects $Ob\,\A = \{Y_1,\dots,Y_m\}$, such that
\begin{equation} \label{eq:directedness}
hom_\A(Y_i,Y_j) = \begin{cases} \text{finite-dimensional over $\K$} & i<j, \\
\K e_i \text{ ($e_i$ is the identity element)}& i = j, \\
0 & i>j. \end{cases}
\end{equation}

A (strictly unital, finite-dimensional, right) $\A$-module $\MM$ consists of a collection of finite-dimensional graded $\K$-vector spaces $\MM(Y_j)$, together with maps
\begin{equation} \label{eq:mu-mm}
\MM(Y_{j_d}) \otimes hom_\A(Y_{j_{d-1}},Y_{j_d}) \otimes \cdots \otimes hom_\A(Y_{j_0},Y_{j_1})
\xrightarrow{\mu_\MM^{d+1}} \MM(Y_{j_0})[1-d]
\end{equation}
for all $d \geq 0$. Strict unitality means that $\mu_\MM^2(\cdot,e_k) = id$, and that all the higher order maps $\mu_\MM^{d+1}$, $d \geq 2$, vanish if one of the last $d$ entries is an identity morphism. The $A_\infty$-module equations are
\begin{equation}
\begin{aligned}
&
0 = \sum_j (-1)^{||a_1|| + \cdots + ||a_j||}  \, \mu_\MM^{j+1}(\mu_\MM^{d-j+1}(m,a_d,\dots,a_{j+1}),\dots, a_1) + \\
& + \sum_{i,j} (-1)^{||a_1|| + \cdots + ||a_j||} \, \mu_\MM^{d-i+2}(m,a_d,\dots,a_{j+i+1},\\[-1em]
& \qquad \qquad \qquad \qquad \qquad \qquad \qquad \mu_\A^i(a_{j+i},\dots,a_{j+1}),
\dots,a_1),
\end{aligned}
\end{equation}
where $||a|| = deg(a) -1$ is the reduced degree. $A_\infty$-modules of this kind form a differential graded category $\CC = mod(\A)$ (note that here, differential graded categories are considered as a special class of $A_\infty$-categories, with the resulting notation and sign conventions). An element $\phi \in hom_\CC^k(\MM,\NN)$ consists of a collection of maps
\begin{equation}
\MM(Y_{j_d}) \otimes hom_\A(Y_{j_{d-1}},Y_{j_d}) \otimes \cdots \otimes hom_\A(Y_{j_0},Y_{j_1}) \stackrel{\phi^{d+1}}{\longrightarrow} \NN(Y_{j_0})[k-d],
\end{equation}
with the property that $\phi^{d+1}$ vanishes if one of the last $d$ entries is an identity morphism (this requirement, together with directedness, ensures that the $hom_\CC$ spaces are always finite-dimensional). The differential is
\begin{equation}
\begin{aligned}
 & \mu^1_\CC(\phi)^{d+1}(m,a_d,\dots,a_1) = \\
 & = \sum_j (-1)^{||a_{j+1}|| + \cdots + ||a_d|| + |m|}\, \mu^{j+1}_\NN(\phi^{d-j+1}(m,a_d,\dots,a_{j+1}),\dots,a_1)
 \\
 & + \sum_j (-1)^{||a_{j+1}|| + \cdots + ||a_d|| + |m|}\,
 \phi^{j+1}(\mu_\MM^{d-j+1}(m,a_d,\dots,a_{j+1}),\dots,a_1)
 \\
 & + \sum_{i,j} (-1)^{||a_{j+1}|| + \cdots + ||a_d|| + |m|}\,
 \phi^{d-i+2}(m,a_d,\dots,a_{j+i+1},
 \\[-1em]
 &
 \qquad\qquad \qquad\qquad\qquad\qquad\qquad\qquad\quad
 \mu_\A^i(a_{j+i},\dots, a_{j+1}),\dots,a_1),
\end{aligned}
\end{equation}
$|m| = deg(m)$ being the ordinary unreduced degree; and the composition is
\begin{equation}
\begin{aligned}
& \mu^2_\CC(\psi,\phi)^{d+1}(m,a_d,\dots,a_1) = \\
& = \sum_j
(-1)^{||a_{j+1}|| + \cdots + ||a_d|| + |m|} \,
\psi^{j+1}(\phi^{d-j+1}(m,a_d,\dots,a_{j+1}), \\[-1em]
& \qquad\qquad\qquad\qquad\qquad\qquad\qquad\qquad\quad\qquad\quad
 a_j,\dots,a_1).
\end{aligned}
\end{equation}

Let $C = H^0(\CC)$ be the underlying cohomological category (the chain homotopy category of $A_\infty$-modules). Let $f \in Hom_C(\MM,\NN)$ be a morphism in that category, and $\phi$ a cocycle representing it. The induced map on $\mu^1$-cohomology,
\begin{equation}
H(\phi^1): \bigoplus_k H^*(\MM(Y_k),\mu^1_\MM) \longrightarrow
\bigoplus_k H^*(\NN(Y_k),\mu^1_\NN)
\end{equation}
depends only on $f$. We say that $f$ is a quasi-isomorphism if $H(\phi^1)$ is an isomorphism.
An important property of $A_\infty$-modules is that every quasi-isomorphism can be inverted in $C$ \cite[Section 4]{keller-notes}.

Submodules and quotient modules of $A_\infty$-modules are defined in the obvious way.  Besides that, we will use a few other constructions. First of all, given $\MM \in Ob\,\CC$ and a finite-dimensional chain complex $(Z,\delta_Z)$ of $\K$-vector spaces, one can form the tensor product $\MM \otimes Z$ \cite[Section 3c]{seidel04}. This is defined by setting $(Z \otimes \MM)(Y_j) = Z \otimes \MM(Y_j)$, with
\begin{equation}
\begin{aligned}
& \mu^1_{Z \otimes \MM}(z \otimes m) = (-1)^{||m||} \delta_Z(z) \otimes m +
z \otimes \mu^1_\MM(m), \\
& \mu^{d+1}_{Z \otimes \MM}(z \otimes m, a_d,\dots,a_1) = z \otimes \mu^{d+1}_\MM(m,a_d,\dots,
a_1) \quad \text{for $d > 0$.}
\end{aligned}
\end{equation}
Note that $Z \otimes \MM$ is quasi-isomorphic to $H^*(Z,\delta_Z) \otimes \MM$, which in turn is a direct sum of shifted copies of $\MM$ (indexed by generators of $H(Z)$, with shifts given by the degress of those generators). Next, given two $A_\infty$-modules and a degree zero cocycle $\phi \in hom^0_\CC(\MM,\NN)$, $\mu^1_\CC(\phi) = 0$, we can form its mapping cone $\QQ = Cone(\phi)$, which is $\QQ(Y_j) = \MM(Y_j)[1] \oplus \NN(Y_j)$, with module structure
\begin{equation}
\begin{aligned}
\mu^{d+1}_\QQ(m \oplus n,a_d,\dots,a_1) & =
\mu^{d+1}_\MM(m,a_d,\dots,a_1) \\ & \!\!\!\!\!\!\! \oplus (\mu^{d+1}_\NN(n,a_d,\dots,a_1) +
\phi^{d+1}(m,a_d,\dots,a_1)).
\end{aligned}
\end{equation}
One can prove that the isomorphism class of $\QQ$ in $C$ depends only on $f = [\phi] \in Hom_C(\MM,\NN)$. Finally, there is a combination of the two last-mentioned operations which will be useful for our purposes. Namely, given $\MM$ and $\NN$, there is a canonical evaluation map $\epsilon \in hom^0_\CC(hom_\CC(\MM,\NN) \otimes \MM, \NN),$ given by
\begin{equation}
\epsilon^{d+1}(\phi \otimes m,a_d,\dots,a_1) =
\phi^{d+1}(m,a_d,\dots,a_1).
\end{equation}
We denote its cone by $T_\MM(\NN) = Cone(\epsilon)$, and call this process (algebraically) twisting $\NN$ by $\MM$.

$\CC$ is a triangulated $A_\infty$-category, with the standard exact triangles involving mapping cones \cite[Section 3h]{seidel04}. Hence, $C$ itself is a triangulated category in the classical sense. In particular, the algebraic twist sits in an exact triangle
\begin{equation} \label{eq:twist-triangle}
\xymatrix{
Hom^*_C(\MM,\NN) \otimes \MM \ar[r] & \NN \ar[r] & T_\MM(\NN)
\ar@/^2pc/[ll]^{[1]}
}
\end{equation}
Here and later on, we write $Hom^i_C(\MM,\NN)$ for the space of degree $i$ morphisms $Hom_C(\MM,\NN[i])$, and $Hom^*_C(\MM,\NN)$ for the direct sum of those spaces over all $i \in \Z$. As a second class of examples, any short exact sequence of modules extends to an exact triangle in $\CC$. This is an analogue of the well-known corresponding property for derived categories of abelian categories, and holds in this context because of the invertibility of quasi-isomorphisms.

\subsection{Simple modules and the canonical filtration}
The smallest nontrivial objects in $\CC$ are the simple modules $\SS_j$, whose underlying vector spaces are
\begin{equation} \label{eq:simple}
\SS_j(Y_i) = \begin{cases} \K \text{ (located in degree $0$)} & i = j, \\
0 & i \neq j. \end{cases}
\end{equation}
In that case, the module structure \eqref{eq:mu-mm} is uniquely determined. Now consider an arbitrary $\MM \in Ob\,\CC$. Because of directedness, this comes with a canonical decreasing filtration, given by the submodules
\begin{equation} \label{eq:m-filtration}
\MM^{\leq m+1-j}(Y_i) = \begin{cases} \MM(Y_i) & i \leq m+1-j, \\ 0 & \text{otherwise.}
\end{cases}
\end{equation}
The graded pieces are $\MM^{\leq m+1-j}/\MM^{\leq m-j} = \MM(Y_{m+1-j}) \otimes \SS_{m+1-j}$. In particular, one gets the following characterization of simple modules.

\begin{Lemma} \label{th:quasi-simple}
Let $\MM$ be an $A_\infty$-module such that all the complexes $(\MM(Y_i),\mu^1_\MM)$ are acyclic except one (say for $i = j$), whose cohomology is one-dimensional and placed in degree zero. Then $\MM \iso \SS_j$ in $C$.
\end{Lemma}

\proof In the canonical filtration, all the graded pieces except one are acyclic, hence isomorphic to zero. From a standard exact triangle argument, it follows that $\MM$ is isomorphic to $\MM^{\leq j}/\MM^{\leq j-1}$, which in turn is isomorphic to $H(\MM(Y_j)) \otimes \SS_j$. \qed

By definition, $hom_\CC(\SS_i,\SS_j) = 0$ whenever $i<j$, and $hom_\CC(\SS_i,\SS_i)$ contains only multiples of the identity map. Passing to the cohomological category $C$, it follows that $(\SS_m,\dots,\SS_1)$ is an exceptional collection. The existence of canonical filtrations shows that this collection is full, which means that the $\SS_j$ generate $C$ as a triangulated category (see \cite{goro} for definitions and further discussion).

\begin{Lemma} \label{th:map-into-simple}
Suppose that $Hom_C^*(\MM,\SS_j)$ is zero for all $j$. Then $\MM$ itself is isomorphic to the zero object in $C$.
\end{Lemma}

\proof This is actually a general property of full exceptional collections. Alternatively, one can argue as follows. Suppose that $\MM$ is nonzero, and take the largest $j$ such that $\MM(Y_j)$ is not acyclic. Looking at the canonical filtration, it follows that $\MM$ is quasi-isomorphic to $\MM^{\leq j}$. Using the assumption, one finds that
\begin{equation}
0 = Hom_C^*(\MM,\SS_j) = Hom_C^*(\MM^{\leq j},\SS_j) = H^*(\MM(Y_j)^\vee),
\end{equation}
which is a contradiction. \qed

\begin{Lemma} \label{th:algebraic-cone}
Let $c \in Hom_C(\MM_1,\MM_0)$ be a morphism such that composition with $c$ yields an isomorphism $Hom_C^*(\MM_0,\SS_j) \rightarrow Hom_C^*(\MM_1,\SS_j)$ for all $j$. Then $c$ itself is an isomorphism.
\end{Lemma}

\proof Choose a cochain representative of $c$, and let $\MM$ be its mapping cone. From the standard exact triangle involving that cone, one gets a long exact sequence
\begin{equation}
 \cdots Hom_C^*(\MM,\SS_j) \rightarrow Hom_C^*(\MM_1,\SS_j) \rightarrow Hom_C^*(\MM_0,\SS_j) \cdots
\end{equation}
where the second $\rightarrow$ is composition with $c$. In view of that, our assumption implies that $Hom_C^*(\MM,\SS_j) = 0$ for all $j$, and then by Lemma \ref{th:map-into-simple} $\MM$ itself is zero. Again appealing to standard facts about mapping cones, it follows that $c$ is an isomorphism. \qed

\subsection{Projective modules and the Yoneda embedding\label{subsec:projective}}
Another basic class of objects in $\CC$ are the elementary projective modules $\PP_k$, given by
\begin{equation}
\PP_k(Y_j) = hom_\A(Y_j,Y_k)
\end{equation}
and $\mu_{\PP_k}^{d+1} = \mu_\A^{d+1}$. For any $\MM$ there is a canonical quasi-isomorphism ${\mathcal F}^1: \MM(Y_k) \rightarrow hom_\CC(\PP_k,\MM)$, given by
\begin{equation} \label{eq:theta}
{\mathcal F^1}(m)^{d+1}(a,a_d,\dots,a_1,a) = \mu_\MM^{d+2}(m,a,a_d,\dots,a_1).
\end{equation}
It is elementary that this is a chain homomorphism. For a proof that it is a quasi-isomorphism, see \cite[\S 7]{fukaya1} or \cite[Section 2g]{seidel04}. An explicit quasi-inverse takes $\phi \in hom_\CC(\PP_k,\MM)$ to $\phi^1(e_k) \in \MM(Y_k)$. Specializing to $\MM = \PP_l$, one has quasi-isomorphisms $hom_\A(Y_k,Y_l) \htp hom_\CC(\PP_k,\PP_l)$ for all $k,l$. This observation can be sharpened as follows. There is a canonical $A_\infty$-functor $\F: \A \rightarrow \CC$ (the Yoneda embedding), sending $Y_k$ to $\PP_k$, which extends the previous $\F^1$, hence is a quasi-equivalence onto its image \cite[\S 9]{fukaya1}.

It follows from the previous discussion that $(\PP_1,\dots,\PP_m)$ is an exceptional collection in the cohomological category $C$. As in the case of simple modules, this collection is full. To prove that, take some $\MM$ with the property that $\MM(Y_j)$ is acyclic for $j > k$, and consider the twisted object $\NN = T_{\PP_k}\MM$. For $j > k$ we have $\PP_k(Y_j) = 0$, hence $\NN(Y_j) = \MM(Y_j)$ remains acyclic. Moreover, $\NN(Y_k)$ is the mapping cone (in the sense of chain complexes) of a map $hom_\CC(\PP_k,\MM) \otimes \PP_k(Y_k) = hom_\CC(\PP_k,\MM) \rightarrow \MM(Y_k)$, which is in fact precisely the inverse quasi-isomorphism described above. The outcome is that $\NN(Y_j)$ is acyclic for $j > k-1$. By repeating this process, one writes an arbitrary module as an iterated mapping cone involving only shifted copies of the $\PP_k$ as building blocks. From this fact and the Yoneda embedding, it then follows that $\CC$ is quasi-equivalent to the derived $A_\infty$-category ${\mathcal D}(\A)$ \cite[Section 5n]{seidel04}. In particular, descending to cohomology, we have an equivalence of triangulated categories
\begin{equation} \label{eq:c-d}
C \iso D(\A^{opp}).
\end{equation}

Finally, we need to discuss briefly the relation between simple and projective modules. This will be based on the theory of mutations in triangulated categories \cite{goro}, which defines an action of the braid group $Br_m$ on the set of full exceptional collections (up to isomorphism) in the category $C$. The standard generators $\sigma_i$, $1 \leq i < m$, act by elementary mutations
\begin{equation}
(\MM_1,\dots,\MM_m) \longmapsto (\MM_1,\dots,\MM_{i-1},T_{\MM_i}(\MM_{i+1}),\MM_i,\MM_{i+2},\dots,\MM_m).
\end{equation}
Suppose that we have two collections $(\MM_1^!,\dots,\MM_m^!)$ and $(\MM_m,\dots,\MM_1)$, of which the second is obtained from the first one through the action of the element $\Delta^{1/2} = \sigma_{m-1}(\sigma_{m-2}\sigma_{m-1}) \cdots (\sigma_1\sigma_2 \cdots \sigma_{m-1}) \in Br_m$. Then these collections are duals \cite[Section 5k]{seidel04} (for earlier work see \cite[Section 2.6]{gorodentsev-kuleshov04} or \cite[Section 7]{bondal89}), in the sense that
\begin{equation} \label{eq:dual-exceptional}
Hom^*_C(\MM_k^!,\MM_j) = \begin{cases} \K \text{ (located in degree $0$)} & j=k, \\
0 & \text{otherwise.}
\end{cases}
\end{equation}
In particular, if $\MM_k^! = \PP_k$ is the projective collection, \eqref{eq:dual-exceptional} determines $H^*(\MM_j(Y_k))$, which in view of Lemma \ref{th:quasi-simple} means that $\MM_j \iso \SS_j$ is isomorphic to the collection consisting of simple modules (there is also a direct, but more computational, proof of this \cite[Section 5o]{seidel04}).

\subsection{The Beilinson spectral sequence}
Given two objects $(\MM_0,\MM_1)$ of $\CC$, the decreasing filtration \eqref{eq:m-filtration} of $\MM = \MM_1$ induces an increasing filtration of $hom_\CC(\MM_0,\MM_1)$, which then gives rise to a spectral sequence converging to $Hom_C^*(\MM_0,\MM_1)$. Using the quasi-isomorphism $\F^1$, we write the starting page of this spectral sequence as
\begin{equation} \label{eq:algebraic-ss}
\begin{aligned}
E_1^{jk} & = Hom^{j+k}_C(\MM_0,\MM_1^{\leq m+1-j}/\MM_1^{\leq m-j}) \\
 & \iso \Big(H^*(\MM_1(Y_{m+1-j})) \otimes Hom_C^*(\MM_0,\SS_{m+1-j})\Big)^{j+k}\\
 & \iso \Big(Hom_C^*(\PP_{m+1-j},\MM_1) \otimes Hom_C^*(\MM_0,\SS_{m+1-j})\Big)^{j+k}.
\end{aligned}
\end{equation}

\begin{Lemma} \label{th:algebraic-edge}
If one identifies $\PP_1 = \SS_1$, the right-sided edge map of \eqref{eq:algebraic-ss} turns into the product $Hom_C^*(\PP_1,\MM_1) \otimes Hom_C^*(\MM_0,\PP_1) \rightarrow Hom_C^*(\MM_0,\MM_1)$.
\end{Lemma}

\proof By definition, the edge homomorphism is induced by the chain map
\begin{equation} \label{eq:chain-edge}
\begin{aligned}
 & hom_\CC(\PP_1,\MM_1) \otimes hom_\CC(\MM_0,\SS_1)
 \\ & = hom_\CC(\MM_0, hom_\CC(\PP_1,\MM_1)  \otimes \SS_1)
 \\ & = hom_\CC(\MM_0, \MM_1(Y_1) \otimes \SS_1)
 \\ &= hom_\CC(\MM_0,\MM_1^{\leq 1})
 \longrightarrow hom_\CC(\MM_0,\MM_1),
\end{aligned}
\end{equation}
where the last step comes from the inclusion $\MM_1^{\leq 1} \rightarrow \MM_1$. Concretely, an element $\psi \in hom_\CC(\PP_1,\MM_1)$ is determined by its first order component, which is a single element $\psi^1 \in \MM_1(Y_1)$, and this is precisely the identification $hom_\CC(\PP_1,\MM_1) = \MM_1(Y_1)$ occurring in \eqref{eq:chain-edge}. In contrast, a $\phi \in hom_\CC(\MM_0,\SS_1)$ consists of a whole series of maps $\phi^{d+1}: \MM_0(Y_{j_d}) \otimes \cdots \otimes hom_\A(Y_{j_0},Y_{j_1}) \rightarrow \K$, for all $j_0 = 1 \leq \cdots \leq j_d$. The map \eqref{eq:chain-edge} takes $\psi \otimes \phi$ to the morphism $\eta$ whose components are
\begin{equation}
\eta^d(m,a_d,\dots,a_1) = \begin{cases}
\psi^1(e_1) \phi^{d+1}(m,a_d,\dots,a_1) & \text{if $j_0 = 1$}, \\
0 & \text{otherwise.}
\end{cases}
\end{equation}
By comparing this with the definition, one sees that after identifying $\PP_1 = \SS_1$, this is indeed just the composition map $\mu^2_\CC$ (to make the signs agree, note that for degree reasons, $\phi^{d+1}(m,a_d,\dots,a_1)$ can only be nonzero if $||a_1|| + \cdots + ||a_d|| + |m| = 0$). \qed

\subsection{The Fukaya category}
Take a Lefschetz fibration $p: X \rightarrow \R$, where the total space is equipped with a complex volume form $\eta$. \cite[Section 18]{seidel04} introduces the Fukaya category $\F(p)$. Objects of this $A_\infty$-category are all closed admissible Lagrangian submanifolds $L \subset X$, as well as the Lefschetz thimbles $\Delta \subset X$ for all vanishing paths $\gamma$ such that for $r \gg 0$, $\gamma'(r) = -i$ (of course, all of them have to be equipped with gradings). Take two dual bases of vanishing paths, and then bend the $\gamma_k^!$ clockwise at infinity until they satisfy the restriction we have just mentioned, while still remaining to the right of the $\gamma_j$ (Figure \ref{fig:bent}). Then, the associated Lefschetz thimbles $\Delta_j$ and $\Delta_k^!$ all become objects of $\F(p)$. Because of isotopy invariance, bending the $\gamma_k^!$ leaves the Floer cohomology groups $HF^*(\tilde\Delta_k^!,\tilde{L})$ and $HF^*(\tilde\Delta_k^!,\tilde\Delta_j)$ unchanged (which is why we choose to keep the notation). These groups, as well as $HF^*(\tilde{L}_0,\tilde{L}_1)$ and $HF^*(\tilde{L},\tilde\Delta_j)$, are the morphisms between the corresponding objects in $H(\F(p))$. This is obvious from the definition in some cases, and proved in \cite[Remark 18.12]{seidel04} in the remaining ones.

It is a nontrivial fact \cite[Propositions 18.17 and 18.23]{seidel04} that $(\tilde\Delta^!_m,\dots,\tilde\Delta^!_1)$ is a full exceptional collection in the derived category $D(\F(p))$. In particular, if $\A \subset \F(p)$ is the directed $A_\infty$-subcategory associated to this collection, there is an $A_\infty$-functor $\A \rightarrow \F(p)$ sending $Y_k$ to $\tilde\Delta_{m+1-k}^!$, and this induces an equivalence of triangulated categories $D(\A) \iso D(\F(p))$ \cite[Theorem 18.24]{seidel04}. In view of \eqref{eq:c-d} we therefore have
\begin{equation} \label{eq:c}
D(\F(p)) \iso C,
\end{equation}
where as usual $C = H^0(\CC)$ is the homotopy category of $A_\infty$-modules over $\A$. By construction, this equivalence sends each $\tilde\Delta_k^!$ to the elementary projective module $\PP_{m+1-k}$. Furthermore, \cite[Proposition 18.23]{seidel04} shows that Hurwitz moves on bases of vanishing paths give rise to mutations of exceptional collections. It is an easy geometric exercise to show that the basis $(\gamma_1,\dots,\gamma_m)$ is obtained from $(\gamma^!_m,\dots,\gamma^!_1)$ by applying the Hurwitz move which corresponds to $\Delta^{1/2} \in Br_m$. From this and the discussion in Section \ref{subsec:projective}, it follows that the equivalence \eqref{eq:c} sends $\tilde\Delta_j$ to $\SS_{m+1-j}$ (strictly speaking, Hurwitz moves do not take gradings of Lefschetz thimbles into account, but it follows from \eqref{eq:dual-convention} that the grading of $\tilde\Delta_j$ is indeed the correct one). In particular, by applying \eqref{eq:algebraic-ss} to the image of two closed admissible Lagrangian submanifolds under \eqref{eq:c}, one immediately obtains the spectral sequence as stated in Theorem \ref{th:ss}. Lemma \ref{th:algebraic-cone} turns into Addendum \ref{th:cone}, and Lemma \ref{th:algebraic-edge} into Addendum \ref{th:edge}. In the latter case, the isotopy between $\Delta_m$ and $\Delta_m^!$ gives rise to an isomorphism between objects in $D(\F(p))$, hence also their images in $C$. It is not a priori clear that this agrees with the purely algebraic identification $\SS_1 = \PP_1$ used in Lemma \ref{th:algebraic-edge}. However, the two can only differ by an element of $\K^\times = Aut_C(\PP_1)$, which is why we pick up a slight ambiguity when translating the result into geometric terms (closer inspection would show that the two isomorphisms actually agree, but we do not need this).

\begin{Remark}
For the benefit of readers wishing to compare the spectral sequence derived here with the formulation in \cite[Corollary 18.27]{seidel04}, we list the differences regarding notation and conventions. First of all, the notion of (exact) Lefschetz fibration in \cite[Section 15d]{seidel04} does not agree exactly with the one here. However, our Lefschetz fibrations can be brought into the form required by \cite{seidel04} through some easy modifications (changing the connection so that parallel transport maps become trivial at infinity, and then restricting to a suitable large compact subset of $X$; compare the discussion in \cite[Section 19d]{seidel04}). Next, the notation for dual bases of vanishing cycles has been swapped. Finally, our ordering of the vanishing paths differs from the one in \cite[Figure 18.19]{seidel04}.
\end{Remark}
\begin{figure}[t]
\begin{centering}
\setlength{\unitlength}{0.00062500in}
\begingroup\makeatletter\ifx\SetFigFont\undefined%
\gdef\SetFigFont#1#2#3#4#5{%
  \reset@font\fontsize{#1}{#2pt}%
  \fontfamily{#3}\fontseries{#4}\fontshape{#5}%
  \selectfont}%
\fi\endgroup%
{\renewcommand{\dashlinestretch}{30}
\begin{picture}(4962,4536)(0,-10)
\put(2450,912){\makebox(0,0)[lb]{{\SetFigFont{10}{12}{\rmdefault}{\mddefault}{\updefault}$\gamma_4$}}}
\drawline(3900.000,2412.000)(3894.944,2514.918)(3879.825,2616.845)
	(3854.787,2716.799)(3820.074,2813.818)(3776.017,2906.967)
	(3723.043,2995.349)(3661.661,3078.113)(3592.462,3154.462)
	(3516.113,3223.661)(3433.349,3285.043)(3344.967,3338.017)
	(3251.818,3382.074)(3154.799,3416.787)(3054.845,3441.825)
	(2952.918,3456.944)(2850.000,3462.000)(2747.082,3456.944)
	(2645.155,3441.825)(2545.201,3416.787)(2448.182,3382.074)
	(2355.033,3338.017)(2266.651,3285.043)(2183.887,3223.661)
	(2107.538,3154.462)(2038.339,3078.113)(1976.957,2995.349)
	(1923.983,2906.967)(1879.926,2813.818)(1845.213,2716.799)
	(1820.175,2616.845)(1805.056,2514.918)(1800.000,2412.000)
\drawline(4350.000,2412.000)(4346.311,2513.395)(4335.795,2614.310)
	(4318.498,2714.286)(4294.499,2812.869)(4263.908,2909.609)
	(4226.864,3004.066)(4183.536,3095.811)(4134.120,3184.426)
	(4078.842,3269.507)(4017.953,3350.668)(3951.731,3427.538)
	(3880.476,3499.768)(3804.513,3567.030)(3724.188,3629.016)
	(3639.866,3685.446)(3551.931,3736.062)(3460.784,3780.634)
	(3366.839,3818.958)(3270.524,3850.862)(3172.276,3876.198)
	(3072.544,3894.853)(2971.782,3906.741)(2870.446,3911.807)
	(2769.000,3910.030)(2667.905,3901.417)(2567.620,3886.007)
	(2468.603,3863.870)(2371.304,3835.108)(2276.165,3799.850)
	(2183.620,3758.258)(2094.090,3710.521)(2007.982,3656.856)
	(1925.689,3597.507)(1847.584,3532.745)(1774.024,3462.864)
	(1705.343,3388.182)(1641.854,3309.039)(1583.845,3225.796)
	(1531.581,3138.831)(1485.300,3048.539)(1445.212,2955.333)
	(1411.499,2859.636)(1384.316,2761.884)(1363.785,2662.521)
	(1350.001,2562.000)
\drawline(4950.000,2412.000)(4947.452,2512.745)(4940.080,2613.253)
	(4927.901,2713.292)(4910.943,2812.633)(4889.244,2911.047)
	(4862.856,3008.308)(4831.837,3104.193)(4796.260,3198.482)
	(4756.207,3290.959)(4711.769,3381.410)(4663.048,3469.628)
	(4610.156,3555.410)(4553.215,3638.559)(4492.356,3718.885)
	(4427.717,3796.203)(4359.449,3870.335)(4287.707,3941.111)
	(4212.657,4008.369)(4134.470,4071.953)(4053.327,4131.719)
	(3969.413,4187.527)(3882.922,4239.252)(3794.052,4286.773)
	(3703.007,4329.981)(3609.996,4368.778)(3515.234,4403.073)
	(3418.937,4432.789)(3321.327,4457.858)(3222.628,4478.220)
	(3123.067,4493.831)(3022.872,4504.653)(2922.273,4510.662)
	(2821.503,4511.845)(2720.791,4508.198)(2620.370,4499.729)
	(2520.470,4486.460)(2421.320,4468.419)(2323.149,4445.648)
	(2226.181,4418.200)(2130.640,4386.137)(2036.744,4349.534)
	(1944.711,4308.474)(1854.750,4263.051)(1767.069,4213.371)
	(1681.869,4159.546)(1599.345,4101.701)(1519.688,4039.969)
	(1443.080,3974.491)(1369.697,3905.418)(1299.708,3832.909)
	(1233.273,3757.129)(1170.546,3678.253)(1111.669,3596.462)
	(1056.779,3511.945)(1006.002,3424.895)(959.453,3335.511)
	(917.241,3244.000)(879.461,3150.572)(846.201,3055.441)
	(817.538,2958.826)(793.536,2860.948)(774.251,2762.033)
	(759.728,2662.307)(750.000,2562.000)
\put(750,2337){\circle*{100}}
\put(750,1287){\circle*{100}}
\put(1800,1662){\circle*{100}}
\put(2400,1662){\circle*{100}}
\drawline(750,1287)(750,837)
\drawline(150,1137)(150,837)
\dashline{60.000}(150,312)(150,12)
\dashline{60.000}(750,312)(750,12)
\drawline(150,837)(150,312)
\drawline(750,837)(750,312)
\dashline{60.000}(1800,312)(1800,12)
\drawline(2400,2412)(2400,312)
\dashline{60.000}(2400,312)(2400,12)
\drawline(1800,312)(1800,2412)
\drawline(3900,2412)(3900,312)
\drawline(4350,2412)(4350,312)
\drawline(4950,2412)(4950,312)
\dashline{60.000}(3300,312)(3300,12)
\dashline{60.000}(3900,312)(3900,12)
\dashline{60.000}(4350,312)(4350,12)
\dashline{60.000}(4950,312)(4950,12)
\drawline(3300,2412)(3300,312)
\drawline(150,1137)(150,1140)(151,1148)
	(152,1160)(154,1178)(157,1199)
	(161,1223)(165,1248)(170,1273)
	(176,1297)(183,1320)(191,1343)
	(200,1365)(211,1388)(223,1412)
	(237,1437)(248,1455)(259,1473)
	(271,1492)(284,1512)(298,1533)
	(313,1555)(328,1578)(345,1601)
	(361,1625)(379,1650)(396,1674)
	(414,1700)(432,1725)(450,1749)
	(468,1774)(485,1798)(502,1821)
	(518,1844)(533,1866)(548,1887)
	(562,1907)(576,1926)(588,1944)
	(600,1962)(615,1984)(628,2006)
	(641,2027)(653,2047)(664,2067)
	(674,2086)(684,2105)(693,2123)
	(700,2141)(707,2158)(713,2175)
	(719,2190)(723,2206)(727,2220)
	(730,2234)(733,2248)(735,2261)
	(737,2274)(740,2290)(741,2306)
	(743,2322)(744,2340)(745,2361)
	(746,2383)(747,2407)(748,2433)
	(748,2461)(749,2488)(749,2513)
	(750,2534)(750,2549)(750,2558)
	(750,2561)(750,2562)
\drawline(750,1287)(750,1290)(751,1296)
	(751,1307)(752,1324)(754,1345)
	(757,1370)(760,1398)(763,1427)
	(767,1457)(772,1486)(777,1514)
	(783,1540)(790,1565)(797,1589)
	(805,1612)(815,1634)(825,1656)
	(837,1677)(850,1699)(862,1718)
	(875,1738)(889,1757)(904,1778)
	(919,1799)(936,1821)(954,1843)
	(972,1866)(991,1889)(1010,1913)
	(1030,1936)(1050,1960)(1070,1984)
	(1090,2008)(1109,2031)(1128,2054)
	(1146,2076)(1164,2098)(1181,2118)
	(1196,2138)(1211,2158)(1225,2177)
	(1238,2195)(1250,2212)(1264,2235)
	(1277,2256)(1288,2278)(1298,2300)
	(1307,2322)(1315,2346)(1321,2370)
	(1327,2396)(1333,2423)(1337,2450)
	(1341,2477)(1344,2501)(1347,2523)
	(1348,2540)(1349,2552)(1350,2559)(1350,2562)
\put(250,1812){\makebox(0,0)[lb]{{\SetFigFont{10}{12}{\rmdefault}{\mddefault}{\updefault}$\gamma_1$}}}
\put(800,912){\makebox(0,0)[lb]{{\SetFigFont{10}{12}{\rmdefault}{\mddefault}{\updefault}$\gamma_2$}}}
\put(800,3462){\makebox(0,0)[lb]{{\SetFigFont{10}{12}{\rmdefault}{\mddefault}{\updefault}$\gamma_1^!$}}}
\put(1575,3537){\makebox(0,0)[lb]{{\SetFigFont{10}{12}{\rmdefault}{\mddefault}{\updefault}$\gamma_2^!$}}}
\put(2300,3057){\makebox(0,0)[lb]{{\SetFigFont{10}{12}{\rmdefault}{\mddefault}{\updefault}$\gamma_3^!$}}}
\put(2465,2187){\makebox(0,0)[lb]{{\SetFigFont{10}{12}{\rmdefault}{\mddefault}{\updefault}$\gamma_4^!$}}}
\put(1850,912){\makebox(0,0)[lb]{{\SetFigFont{10}{12}{\rmdefault}{\mddefault}{\updefault}$\gamma_3$}}}
\drawline(3300.000,2412.000)(3288.718,2512.134)(3255.436,2607.248)
	(3201.824,2692.570)(3130.570,2763.824)(3045.248,2817.436)
	(2950.134,2850.718)(2850.000,2862.000)(2749.866,2850.718)
	(2654.752,2817.436)(2569.430,2763.824)(2498.176,2692.570)
	(2444.564,2607.248)(2411.282,2512.134)(2400.000,2412.000)
\end{picture}
}
\caption{\label{fig:bent}}
\end{centering}
\end{figure}

\end{document}